\documentclass{aptpub}

\usepackage[colorlinks=true, linkcolor=blue, urlcolor=blue, citecolor=blue, anchorcolor=blue]{hyperref}
\usepackage{mathtools}
\usepackage{pdfrender}
\usepackage{enumitem}
\usepackage{listings}
\usepackage{xspace}

\usepackage{xypic}
\usepackage[all]{xy}
\xyoption{all}
\xyoption{2cell}
\UseTwocells
\xyoption{tips}
\SelectTips{cm}{}  

\usepackage{tikzit}
\usetikzlibrary{external}
\newcommand{\ignore}[1]{}

\newif\ifignore 
\ignorefalse

\newcommand{\auxproof}[1]{
\ifignore\mbox{}\newline
\textbf{PROOF:} \dotfill\newline
{\it #1}\mbox{}\newline
\textbf{ENDPROOF}\dotfill
\fi}

\newenvironment{myproof}[1][Proof]%
   { \begin{trivlist}%
     \item[\hskip \labelsep {\it #1}.]%
   }%
   { \end{trivlist}%
   }

\newcommand*{\fatten}[1][.4pt]{%
  \textpdfrender{
    TextRenderingMode=FillStroke,
    LineWidth={\dimexpr(#1)\relax},
  }%
}

\ifdefined\mathsl\else
  \DeclareMathAlphabet{\mathsl}{\encodingdefault}{\rmdefault}{\mddefault}{\sldefault}
  \SetMathAlphabet{\mathsl}{bold}{\encodingdefault}{\rmdefault}{\bfdefault}{\sldefault}
\fi

\makeatletter
\newcommand{\mathoverlap}[2]{\mathpalette\mathoverlap@{{#1}{#2}}}
\newcommand{\mathoverlap@}[2]{\mathoverlap@@{#1}#2}
\newcommand{\mathoverlap@@}[3]{\ooalign{$\m@th#1#2$\crcr\hidewidth$\m@th#1#3$\hidewidth}}
\makeatother

\newcommand{\NNO}{\mathbb{N}}

\newcommand{\pR}{\mathbb{R}_{>0}}
\newcommand{\nnR}{\mathbb{R}_{\geq 0}}

\newcommand{\finset}[1]{\ensuremath{\boldsymbol{#1}}}

\newcommand{\setin}[3]{\{#1\in#2\;|\;#3\}}

\newcommand{\bigsetin}[3]{\big\{#1\in#2\;\big|\;#3\big\}}
\newcommand{\Bigset}[2]{\Big\{{\kern.2em}#1\;\Big|\;#2{\kern.2em}\Big\}}
\newcommand{\Bigsetin}[3]{\Big\{{\kern.2em}#1\in#2\;\Big|\;#3{\kern.2em}\Big\}}

\newcommand{\setsize}[1]{|{\kern.1em}#1{\kern.1em}|}
\newcommand{\bigsetsize}[1]{\big|{\kern.1em}#1{\kern.1em}\big|}
\newcommand{\Bigsetsize}[1]{\Big|{\kern.1em}#1{\kern.1em}\Big|}
\newcommand{\one}{\ensuremath{\mathbf{1}}}
\newcommand{\zero}{\ensuremath{\mathbf{0}}}

\newcommand{\after}{\mathrel{\circ}}
\newcommand{\klafter}{\mathbin{\mathoverlap{\circ}{\cdot}}}
\newcommand{\idmap}[1][]{\ensuremath{\mathsl{id}_{#1}}}

\newcommand{\ket}[1]{\ensuremath{|{\kern.1em}#1{\kern.1em}\rangle}}
\newcommand{\bigket}[1]{\ensuremath{\big|{\kern.1em}#1{\kern.1em}\big\rangle}}
\newcommand{\ketstrut}{\vrule height 8.5pt depth 4.5pt width 0pt}
\newcommand{\Bigket}[1]{\ensuremath{\left|\ketstrut{\kern.1em}\right.{\kern-.2em}#1{\kern-.2em}\left.\ketstrut{\kern0em}\right>}}
\newcommand{\bra}[1]{\ensuremath{\langle{\kern.1em}#1{\kern.1em}|}}

\newcommand{\concat}{\ensuremath{\mathbin{+{\kern-.5ex}+}}}

\newcommand{\support}{\ensuremath{\mathsl{supp}}}
\newcommand{\revert}{\ensuremath{\mathsl{rev}}}
\newcommand{\Revert}{\ensuremath{\mathsl{Rev}}}
\newcommand{\Prob}{\ensuremath{\mathsl{P}}}

\newcommand{\shift}{\ensuremath{\mathsl{shift}}}

\newcommand{\coefm}[1]{\ensuremath{\fatten[0.6pt]{(}{\kern1pt}#1{\kern1pt}\fatten[0.6pt]{)}}}

\newcommand{\bibinom}[2]{\left(\!\binom{#1}{#2}\!\right)}

\newcommand{\acc}{\ensuremath{\mathsl{acc}}}

\newcommand{\flrn}{\ensuremath{\mathsl{Flrn}}}
\newcommand{\unif}{\ensuremath{\mathsl{uni}}}
\newcommand{\plusje}{\ensuremath{{\kern-2pt}+{\kern-2pt}}}
\newcommand{\minnetje}{\ensuremath{{\kern-1.5pt}-{\kern-1.5pt}}}

\newcommand{\multisetcoefficientdistribution}{\ensuremath{\mathsl{mc}}}

\newcommand{\boltzmannmlt}{\ensuremath{\mathsl{BoM}}}
\newcommand{\boltzmannmltmlt}{\ensuremath{\mathsl{BoMM}}}
\newcommand{\boltzmannnum}{\ensuremath{\mathsl{BoN}}}
\newcommand{\boltzmannene}{\ensuremath{\mathsl{BoE}}}

\newcommand{\hypgeom}{\ensuremath{\mathsl{hg}}}
\newcommand{\hypergeometric}[1][]{\ensuremath{\hypgeom[#1]}}
\newcommand{\pol}{\ensuremath{\mathsl{pl}}}
\newcommand{\polya}[1][]{\ensuremath{\pol[#1]}}
\newcommand{\nom}{\ensuremath{\mathsl{nom}}}
\newcommand{\nomial}[1][]{\ensuremath{\nom[#1]}}
\newcommand{\mean}{\ensuremath{\mathsl{Mean}}}
\newcommand{\som}{\ensuremath{\mathsl{sum}}}

\newcommand{\Var}{\ensuremath{\mathsl{Var}}}

\newcommand{\intd}{{\kern.2em}\mathrm{d}{\kern.03em}}

\newcommand{\Mlt}{\ensuremath{\mathcal{M}}}
\newcommand{\fullMlt}{\ensuremath{\Mlt_{{\kern-0.2pt}\mathsl{fs}}}}
\newcommand{\natMlt}{\ensuremath{\mathcal{M}}}
\newcommand{\fullnatMlt}{\ensuremath{\natMlt_{{\kern-0.8pt}\mathsl{fs}}}}
\newcommand{\Dst}{\ensuremath{\mathcal{D}}}

\newcommand{\congrightarrow}{\mathrel{\smash{\stackrel{
           \raisebox{.5ex}{$\scriptstyle\cong$}}{
           \raisebox{0ex}[0ex][0ex]{$\rightarrow$}}}}}

\DeclareFixedFont{\ttb}{T1}{txtt}{bx}{n}{11} 
\DeclareFixedFont{\ttm}{T1}{txtt}{m}{n}{11}  
\usepackage{color}
\definecolor{deepblue}{rgb}{0,0,0.5}
\definecolor{deepred}{rgb}{0.6,0,0}
\definecolor{deepgreen}{rgb}{0,0.5,0}
\definecolor{lightgray}{rgb}{0.83,0.83,0.83}

\newcommand{\Python}{\textrm{Python}\xspace}
\newcommand\pythonstyle{\lstset{
language=Python,
basicstyle=\small\ttfamily,
otherkeywords={self,>>>,...},             
keywordstyle=\ttb\color{deepblue},
emph={MyClass,__init__},          
emphstyle=\ttb\color{deepred},    
stringstyle=\color{deepgreen},
frame=tb,                         
showstringspaces=false            %
}}

\lstnewenvironment{python}[1][]
{
\pythonstyle
\lstset{#1}
}
{}

\newcommand\pythoninline[1]{{\pythonstyle\lstinline!#1!}}
\newcommand\inline[1]{{\lstinline!#1!}}

\authornames{B.~Jacobs} 

\shorttitle{Discrete Boltzmann distributions} 

\begin{document}

\title{Discrete Boltzmann distributions via multisets and their coefficients}

\authorone[Radboud University]{Bart Jacobs} 



\addressone{Institute for Computing and Information Science, Radboud University Nijmegen, The Netherlands} 

\emailone{bart@cs.ru.nl} 

\begin{abstract}
This paper investigates the combinatorics that gives rise to the
Boltzmann probability distribution. Despite being one of the most
important distributions in physics and other fields of science, the
mathematics of the underlying model of particles at different energy
levels is underexplored. This paper gives a reconstruction, using
multisets with fixed sums as mathematical representations. Counting
(the coefficients of) such multisets gives a general description of
binomial, trinomial, quadrinomial etc.\ coefficients, here called
$N$-nomials. These coefficients give rise to multiple discrete
Boltzmann distributions that are linked to explanations in the physics
literature.
\end{abstract}

\keywords{Boltzmann distribution; particles; energy; multiset; coefficient}


\ams{05A10}{60C05; 70A05}    





\section{Introduction}



Ludwig Boltzmann introduced the distribution later named after him in
the late 19th century in his work on statistical mechanics. This
`Boltzmann' distribution describes the energy distribution of atoms
and moleculues. The discrete version of the Boltzmann distribution
describes the probabilities of particle configurations over various
energy levels in an equilibrium, for a fixed total energy, see
\textit{e.g.}~\cite{DillB10} for background information on the usage
in natural sciences, or~\cite{DragulescuY00} for usage in economics.
There are several variations of this Boltzmann distribution, in
discrete and continuous form. The continuous version is an exponential
distribution that is obtained via a limit process, using an
approximation. In general, physicists quickly move to this continuous
formulation, based on a sketch of the underlying discrete
situation. This paper, however, concentrates on this discrete stage
and reconstructs its rich mathematical structure in terms of
`$N$-nomial coefficients' and multisets.

These $N$-nomial coefficients are introduced in this paper as
generalisations of binomial, trinomial and quadrinomial coefficients.
This topic exists in the literature only in sketchy
form. Section~\ref{NomialSeqSec} introduces the $N$-nomial
$C_{N}(K,i)$ in a new way as the number of sequences of length $K$,
containing numbers $n_{1}, \ldots n_{K} \in \{0,1,\ldots,N\minnetje
1\}$, satisfying $n_{1} + \cdots + n_{K} = i$.  The standard binomial
coefficients $\binom{K}{i}$ arise as special cases, for
$N=2$. Trinomial and quadrinomial coefficients arise for $N=3$ and
$N=4$, respectively. In physical terms, $C_{N}(K,i)$ is the number of
`microstates' of size $K$ with total energy~$i$, at energy levels
$\{0,1,\ldots,N\minnetje 1\}$. Several other formulations of
$N$-nomials will be given, including one in algorithmic form.  Also,
it will be shown how these $N$-nomials arise in certain polynomial
expressions --- see Theorem~\ref{PolynomialThm} --- corresponding to
the definition given by the Online Encyclopedia of Integer Sequences
(OEIS).

The paper makes systematic use of multisets. These can be seen as
`sets' in which elements may occur multiple times, or as `lists' in
which the order of the elements does not matter. These multisets are
used for two purposes, namely for:
\begin{itemize}
\item an alternative description of $N$-nomials, see
  Section~\ref{NomialMltSec};

\item capturing multiple particles at different energy levels.
\end{itemize}

\noindent Multisets are not explicitly used in statistical mechanics,
but they occur implicitly and can be recognised as such. For instance,
they are called compositions in~\cite[Chapt.~2]{DillB10}.  There, the
probability of 4 tosses of a fair coin are discussed for sequences
$HHHH$ and $HTHH$, where $H$ and $T$ stand for heads and tails. The
probability is $\big(\frac{1}{2}\big)^{4} = \frac{1}{16}$ for both
sequences. Then: ``Contrast this with a question of composition. Which
\emph{composition} is more probable: 4~$H$'s and 0~$T$'s, or 3~$H$'s
and 1~$T$?'' Using the notation introduced in
Section~\ref{NomialMltSec} we will write these compositions as two
multisets $4\ket{H}$ and $3\ket{H} + 1\ket{T}$.  They have
(multinomial) probabilities $\frac{1}{16}$ and $\frac{1}{4}$,
respectively. Some authors also use `microstates' for sequences and
`macrostates' for multisets, but the latter notion is not explicitly
defined.  To quote from~\cite[\S1.6]{PathriaB11}: ``Thus, the correct
way of specifying a microstate of the system is through the
distribution numbers $\{n_{j}\}$, and not through the statement as to
"which particle is in which state."'' Section~\ref{NomialMltSec}
describes an accumulation operation $\acc$ that turns a sequence into
a multiset, basically by ignoring the order of the elements and only
counting their multiplicities (``distribution numbers''). Using the
above coins, we have $\acc(HHHH) = 4\ket{H}$ and $\acc(HTHH) =
3\ket{H} + 1\ket{T}$.  The probability $\frac{1}{4}$ for the multiset
$3\ket{H} + 1\ket{T}$ is obtained as the sum $\frac{1}{4} =
\frac{1}{16} + \frac{1}{16} + \frac{1}{16} + \frac{1}{16}$ of the four
probabilities of the four sequences $THHH$, $HTHH$, $HHTH$, $HHHT$
that accumlate toe $3\ket{H} + 1\ket{T}$. 

\ignore{

print( Multinomial(4)(Flip(Frac(1,2))) )

# 1/16|4|H>> + 
# 1/4|3|H> + 1|T>> + 
# 3/8|2|H> + 2|T>> + 
# 1/4|1|H> + 3|T>> + 
# 1/16|4|T>>

}

The paper introduces multiple versions of a discrete
`Boltzmann' distribution.
\begin{enumerate}
\item A \emph{Boltzmann-on-multisets} distribution appears in
  Section~\ref{BoltzmannMltSec} as a distribution over multisets. It
  captures the distribution of configurations of a fixed number of
  particles over different energy levels, with a fixed total
  energy. It is shown that this distribution corresponds directly to
  exemplaric explations in the physics literature. The Boltzmann
  distribution on multisets (macrostates) arises as image distribution
  (via accumulation) from the uniform distribution on sequences
  (microstates).

\item A \emph{Boltzmann-on-numbers} distribution is obtained from this
  Boltzmann-on-multi\-sets distribution, in
  Section~\ref{BoltzmannNumSec}, as a distribution over the various
  energy levels, again for a fixed total energy. We use `frequentist
  learning' as the canonical operation that turns a multiset into a
  distribution, basically by normalisation (\textit{i.e.}~by
  counting). The Boltzmann-on-numbers distribution arises by applying
  frequentist learning `in probability' to the Boltzmann-on-multisets
  distribution. Several alternative formulations are given for this
  distribution, one of which is entirely in terms of $N$-nomials.

\item Finally, a \emph{Boltzmann-on-energy} distribution is introduced
  as a special case of the Boltzmann-on-numbers distribution, namely
  where the number of energy levels exceeds the total energy. This is
  common in physics. It is shown that this special cases is
  mathematically also very relevant, since it gives a simple
  description of the $N$-nomial, in terms of a multichoose
  coefficient, see Theorem~\ref{BelowThm}~\eqref{BelowThmNomial}.
  This leads to a simple description of this third Boltzmann-on-energy
  distribution that matches formulations in the (physics) literature.
\end{enumerate}

This paper contains several contributions. First, there is the
systematic development of $N$-nomial coefficients, including the
simple description in the physically relevant special case. Next, it
distinguishes the above three discrete Boltzmann distributions, with
their relationship and properties --- such as their symmetry, via
closure under reversal.  Thirdly, the paper establishes a close, new
connection between $N$-nomial coefficients and these three Boltzmann
distributions.  Fourthly, the Boltzmann-on-multisets and
Boltzmann-on-numbers distributions are described as equilibria (fixed
points, stationary distributions) of particular Markov chains, see
Subsection~\ref{MarkovChainSubsec}.  Finally, the paper introduces two
new multivariate distributions via $N$-nomials, exploiting the
Vandermonde property that is proven for $N$-nomials, see
Subsection~\ref{NomialDstSubsec}. One distribution arises similarly to
the multivariate hypergeometric and P\'olya distributions, but differs
from both of them. The other one is a multivariate Boltzmann
distribution that could be used to capture combinations of particles
of different kinds.

This paper builds on ideas and techniques from category theory (see
\textit{e.g.}~\cite{MacLane71,Leinster14}), but mostly hides this
underlying approach, in order not to limit the audience.  The one
place where category theory pops up is in the form of `functoriality'
of the multiset and distribution operations $\natMlt$ and $\Dst$, but
there the relevant notions are introduced concretely. Thus, no
familiarity with category theory is assumed, but readers who are
interested in this wider more abstract perspective are referred
to~\cite{Jacobs21g}.

\section{Multichoose coefficients}\label{MultiChooseSec}

We shall frequently use the \emph{binomial} coefficient and
\emph{multichoose} coefficient, given in terms of factorials
respectively as:
\[ \begin{array}{rclcrcl}
\displaystyle\binom{n}{i}
& = &
\displaystyle\frac{n!}{(n-i)!\cdot i}
& \qquad\qquad &
\displaystyle\bibinom{m}{j}
& = &
\displaystyle\frac{(m+j-1)!}{(m-1)!\cdot j!}.
\end{array} \]

\noindent The constraints are that $n \geq 0$ and $0 \leq i \leq n$
and that $m\geq 1$. It is well-known that $\binom{n}{i}$ is the number
of \emph{subsets} of size $i$ of a set of size $n$.  It is less
well-known that $\big(\binom{m}{j}\big)$ is the number of
\emph{multisets} of size $j$ of a set of size $m$, see
Lemma~\ref{MltLem}~\eqref{MltLemSize} below.

We shall need the following basic equations. The first one can be
found for instance in~\cite[Identity~149]{BenjaminQ03}.  The next two
can be obtained easily from this first equation.

\begin{lemma}
\label{ChooseSumLem}
Let $n\geq 1$.
\begin{enumerate}
\item \label{ChooseSumLemZero} For $m\geq 1$,
\[ \begin{array}{rcl}
\displaystyle\sum_{0\leq j < m}\, \bibinom{n}{j}
& = &
\displaystyle\bibinom{m}{n}.
\end{array} \]

\auxproof{
\cite[Identity~149]{BenjaminQ03} allows $n=0$, without further
explanation.
\[ \begin{array}{rcl}
\displaystyle\sum_{k=0}^{m} \bibinom{n}{k}
& = &
\displaystyle\bibinom{n+1}{m}
\end{array} \]

\noindent Our formulation yields
\[ \begin{array}{rcccccl}
\displaystyle\sum_{k=0}^{m} \bibinom{n}{k}
& = &
\displaystyle\bibinom{m+1}{n}
& = &
\displaystyle\frac{(n+m)!}{m!\cdot n!}
& = &
\displaystyle\bibinom{n+1}{m}
\end{array} \]
}

\item \label{ChooseSumLemOne} For $m\geq 2$,
\[ \begin{array}{rcl}
\displaystyle\sum_{0\leq j < m}\, \bibinom{n}{j} \cdot j
& = &
\displaystyle n\cdot \bibinom{m-1}{n+1}
\end{array} \]

\item \label{ChooseSumLemTwo} For $m\geq 3$,
\[ \begin{array}{rcl}
\displaystyle\sum_{0\leq j < m}\, \bibinom{n}{j} \cdot j^{2}
& = &
\displaystyle n\cdot (n+1) \cdot \bibinom{m-2}{n+2} + n\cdot\bibinom{m-1}{n+1}.
\end{array} \eqno{\square} \]
\end{enumerate}
\end{lemma}

\section{$N$-nomial coefficients via sequences}\label{NomialSeqSec}

This section introduces $N$-nomial coefficients, as generalisations of
binomial, trinomial, quadrinomial \emph{etc.}\ coefficients. There
does not seem to be a systematic description of these coefficients in
the literature. The website~\cite{OEISNnomials} of the Online
Encyclopedia of Integer Sequences (OEIS) contains a description of how
$N$-nomials arise in expansions of certain polynomials, together with
multiple examples for trinomial and quadrinomial coefficients. This
section introduces a new direct formula, for the $N$-ary case. It will
be related to the polynomial formulation below, in
Theorem~\ref{PolynomialThm}. The definition is first given in terms of
sequences, and then reformulated in terms of multisets, in the next
section.  Several properties are included for these $N$-nomials.

We recall that the binomial coefficient $\binom{N}{i} \in \NNO$ gives
the number of subsets of size~$i$, of a (finite) set of size~$N$. Such
a subset of size $i$, say of a set $\{x_{1}, \ldots, x_{N}\}$ of $N$
elements, can be identified with sequence of binary numbers $(b_{1},
\ldots, b_{N}) \in \{0,1\}^{N}$ of length $N$ with numbers $b_{i} \in
\{0,1\}$ satisfying $i = \som\big(b_{1}, \ldots, b_{N}\big) = b_{1} +
\cdots + b_{N}$.  Each number $b_{i}$ then tells if the element
$x_{i}$ is in the subset (when $b_{i} = 1$) or not (when $b_{i} = 0$).

This description of the binomial coefficient can be generalised to
(univariate) trinomial, quadrinomial \ldots coefficients.  We use
$\setsize{X} \in \NNO$ to describe the size, that is, the number of
elements, of a finite set $X$.

\begin{definition}
\label{NnomialDef}
We fix a natural number $N\geq 1$ and write $\finset{N} =
\{0,1,\ldots,N\minnetje 1\} \subseteq \NNO$ for the (sub)set of the
first $N$ natural numbers. For numbers $K\in\NNO$ and $0 \leq i \leq
(N\minnetje 1)\cdot K$ we define the \textit{$N$-nomial coefficient}
$C_{N}(K,i) \in \NNO$ as:
\begin{equation}
\label{NnomialSeqEqn}
\begin{array}{rcl}
C_{N}(K, i)
& \coloneqq &
\Bigsetsize{\bigsetin{\vec{n}}{\finset{N}^{K}}{\som(\vec{n}) = i}}.
\end{array}
\end{equation}
\end{definition}

From this definition it is obvious that:
\[ \begin{array}{rclcrclcrcl}
C_{N}(0, 0)
& = &
1
& \quad &
C_{1}(K, 0)
& = &
1
& \quad &
C_{2}(K, i)
& = &
\displaystyle\binom{K}{i}, \;\mbox{ for } 0\leq i\leq K.
\end{array} \]

\noindent For $N=3$, definition~\eqref{NnomialSeqEqn} gives the
trinomial coefficients, written as triangle:
\begin{center}
\begin{tabular}{c|ccccccccccc}
$\; K \;$ & & & & & & \makebox[0pt]{$C_{3}(K, -)$} & & & & &
\\
\hline
\hline
0 & & & & & & 1 & & & & & 
\\
1 & & & & & 1 & 1 & 1 & & & & 
\\
2 & & & & 1 & 2 & 3 & 2 & 1 & & &
\\
3 & & & 1 & 3 & 6 & 7 & 6 & 3 & 1 & &
\\
4 & & 1 & 4 & 10 & 16 & 19 & 16 & 10 & 4 & 1 &
\\
$\vdots$ 
\end{tabular}
\end{center}

\noindent The last line says, for instance $C_{3}(4,0) = 1$, 
$C_{3}(4,1) = 4$, $C_{3}(4,2) = 10$, \textit{etc}.

In a similar way, the quadrinomial coefficients are obtained
from~\eqref{NnomialSeqEqn}, now written as sequences:
\[ \begin{array}{rcl}
C_{4}(0, -)
& = &
1
\\
C_{4}(1, -)
& = &
1, 1, 1, 1
\\
C_{4}(2, -)
& = &
1, 2, 3, 4, 3, 2, 1
\\
C_{4}(3, -)
& = &
1, 3, 6, 10, 12, 12, 10, 6, 3, 1
\\
C_{4}(4, -)
& = &
1, 4, 10, 20, 31, 40, 44, 40, 31, 20, 10, 4, 1
\\
C_{4}(5, -)
& = &
1, 5, 15, 35, 65, 101, 135, 155, 155, 135, 101, 65, 35, 15, 5, 1
\\
\vdots\qquad & &
\end{array} \]

\noindent The occurrence of $N$-nomials in certain polynomial
expressions will be described in Theorem~\ref{PolynomialThm}. We first
collect some basic properties.

\begin{lemma}
\label{NnomialSeqLem}
Let numbers $N \geq 1$, $K\geq 0$ be given.
\begin{enumerate}
\item \label{NnomialSeqLemSum} $N$-nomials add up to $N^K$ in the following way.
\[ \begin{array}{rcl}
\displaystyle\sum_{0 \leq i \leq (N-1)\cdot K} C_{N}(K, i)
& = &
N^{K}.
\end{array} \]

\noindent This generalises the familiar equation $\sum_{0\leq i\leq K}
\binom{K}{i} = 2^{K}$ in the binary case.

\item \label{NnomialSeqLemVDM} $N$-nomials satisfy a Vandermonde
  property: for each $0 \leq i \leq (N\minnetje 1)\cdot K$, if $K =
  K_{1} + K_{2}$, then:
\[ \begin{array}{rcl}
C_{N}(K,i)
& = &
\displaystyle\sum_{0 \leq i_{1} \leq (N-1)\cdot K_{1}, \,
   0 \leq i_{2} \leq (N-1)\cdot K_{2}, \, i_{1} + i_{2} = i}
   C_{N}(K_{1}, i_{1})\cdot C_{N}(K_{2}, i_{2}).
\end{array} \]

\noindent The (familiar) Vandermonde property for binomial
coefficients is a special case, for $N=2$.
\end{enumerate}
\end{lemma}

\begin{myproof}
\begin{enumerate}
\item Since:
\[ \begin{array}{rcl}
\displaystyle\sum_{0 \leq i \leq (N-1)\cdot K} \, C_{N}(K, i)
& = &
\displaystyle\sum_{0 \leq i \leq (N-1)\cdot K} \,
   \Bigsetsize{\bigsetin{\vec{n}}{\finset{N}^{K}}{\som(\vec{n}) = i}}
\\[+1.0em]
& = &
\Bigsetsize{\bigsetin{\vec{n}}{\finset{N}^{K}}
   {0 \leq \som(\vec{n}) \leq (N-1)\cdot K}}
\\[+0.4em]
& = &
\Bigsetsize{\finset{N}^{K}} 
\hspace*{\arraycolsep}=\hspace*{\arraycolsep}
N^{K}.
\end{array} \]

\item We use that sequences of length $K = K_{1} + K_{2}$ can be written as
concatenations $\concat$ of sequences of length $K_{1}$ and $K_{2}$.
\[ \hspace*{-1em}\begin{array}[b]{rcl}
\lefteqn{C_{N}(K,i)}
\\
& \smash{\stackrel{\eqref{NnomialSeqEqn}}{=}} &
\Bigsetsize{\bigsetin{\vec{n}}{\finset{N}^{K}}{\som(\vec{n}) = i}}
\\[+0.2em]
& = &
\displaystyle\sum_{0 \leq i_{1} \leq (N-1)\cdot K_{1}, \,
   0 \leq i_{2} \leq (N-1)\cdot K_{2}, \, i_{1} + i_{2} = i}
   \Big|\,\big\{\vec{n_1} \concat \vec{n_2} \; \big| \;
   \vec{n_{1}} \in \finset{N}^{K_1}, \vec{n_{2}} \in \rlap{$\finset{N}^{K_2},$}
\\[-0.8em]
& & \hspace*{17em}
   \som(\vec{n_1}) = i_{1}, \som(\vec{n_2}) = i_{2}\big\}\,\Big|
\\[+0.4em]
& = &
\displaystyle\sum_{0 \leq i_{1} \leq (N-1)\cdot K_{1}, \,
   0 \leq i_{2} \leq (N-1)\cdot K_{2}, \, i_{1} + i_{2} = i}
   \Bigsetsize{\bigsetin{\vec{n_1}}{\finset{N}^{K_1}}{\som(\vec{n_1}) = i_1}}
\\[-0.8em]
& & \hspace*{16.5em}
   \cdot\,
   \Bigsetsize{\bigsetin{\vec{n_2}}{\finset{N}^{K_2}}{\som(\vec{n_2}) = i_2}}
\\[+0.4em]
& = &
\displaystyle\sum_{0 \leq i_{1} \leq (N-1)\cdot K_{1}, \,
   0 \leq i_{2} \leq (N-1)\cdot K_{2}, \, i_{1} + i_{2} = i}
   C_{N}(K_{1}, i_{1})\cdot C_{N}(K_{2}, i_{2}).
\end{array} \eqno{\square} \]
\end{enumerate}
\end{myproof}

\section{$N$-nomial coefficients via multisets}\label{NomialMltSec}

This section introduces a second formulation of $N$-nomials in terms
of multisets, and also a short algorithm for computing $N$-nomials.
The multiset reformulation of $N$-nomials allows us to prove
additional basic properties about them, and in particular, to relate
our formulation of $N$-nomials to their description in terms of
polynomials (see Theorem~\ref{PolynomialThm}). Since multisets are not
so common yet, we start with a brief introduction.

Informally, multisets are `subsets' in which elements may occur
multiple times.  One can also view multisets as `sequences' in which
the order of the elements does not matter. Multisets play a
fundamental role in combinatorics and probability, but this is not
always recognised.  One reason may be that there is no commonly
accepted, good notation for multisets.  We shall use `ket' notation,
borrowed from quantum theory. An urn $\upsilon$ containing three red
($R$), two green ($G$) and one blue ($B$) ball will be written as:
\begin{equation}
\label{Urn}
\begin{array}{rcl}
\upsilon
& = &
3\bigket{R} + 2\bigket{G} + 1\bigket{B}.
\end{array}
\end{equation}

\noindent Some authors, for instance in~\cite{GrahamKP94,BenjaminQ03},
simply repeat elements in subset notation, as in
$\{R,R,R,G,G,B\}$. This does not really scale.

The numbers $3,2,1$ before the three kets in~\eqref{Urn} describe the
multiplicities of the elements $R,G,B$ within the kets
$\ket{-}$. These kets have no mathematical meaning and are only used
to separate the multiplicities and their elements. The convention is
that the order of the kets does not matter, that multiplicities for
the same elements can be added, and that zero-multiplicities can be
omitted. Thus we could also write $\upsilon = 1\ket{R} + 1\ket{B} +
2\ket{G} + 2\ket{R} + 0\ket{P}$, where $P$ stands for purple, say.

Here is another example: the polynomial $x^{3} - 7x^{2} + 16x -12 =
(x-2)(x-2)(x-3)$ has a multiset of roots, namely:
$2\ket{2} + 1\ket{3}$.

In a slightly more formal way, we shall write a multiset with elements
from an arbitrary set $X$ as a function $\varphi \colon X \rightarrow
\NNO$ with finite support: the set $\support(\varphi) \coloneqq
\setin{x}{X}{\varphi(x) \neq 0}$ is required to be finite. Thus we can
write $\varphi = \sum_{x} \varphi(x)\bigket{x}$. We shall freely
switch between the ket notation and this functional notation, and use
whichever is most convenient. Notice that in this context, a multiset
is always finite, but it may be defined over an arbitrary,
not-necessarily finite set $X$.

The \emph{size} of a multiset is the total number of its elements,
including multiplicities. Thus, the urn $\upsilon$ in~\eqref{Urn} has
size~$6$.  In general, we write $\|\varphi\| \coloneqq \sum_{x}
\varphi(x)$ for the size of a multiset $\varphi$. We shall also write
$\natMlt(X)$ for the set of all multisets over $X$, and $\natMlt[K](X)
\subseteq \natMlt(X)$ for the subset of multisets of size
$K\in\NNO$. For $K=0$, the set $\natMlt[K](X)$ has precisely one
member, namely the empty multiset $\zero$ with zero elements. For the
above urn~\eqref{Urn} we can now write $\upsilon \in
\natMlt[6]\big(\{R,G,B\}\big)$.

An easy mental picture for multisets in $\natMlt(X)$ is to think of
them as urns filled with multiple balls with different colours, from
the set $X$.  Putting the contents of two urns together corresponds to
(pointwise) addition of multisets. This makes $\natMlt(X)$ a
commutative monoid, with the empty multiset $\zero$ as neutral
element. In fact, $\natMlt(X)$ is the \emph{free} commutative monoid
on the set $X$.  Another common, but less intuitive, model of
multisets involves `stars and bars', see
\textit{e.g.}~\cite[\S5.4]{BenjaminQ03}.

\begin{definition}
\label{SumDef}
We often consider multisets $\varphi\in\natMlt(\finset{N})$ on the set
$\finset{N} = \{0,1,\ldots,N\minnetje 1\}$ containing the first $N$
natural numbers. For such a multiset
$\varphi\in\natMlt[K](\finset{N})$ we define $\som(\varphi) \in \NNO$
as:
\begin{equation}
\label{SumEqn}
\begin{array}{rcl}
\som(\varphi) 
& \coloneqq &
\displaystyle \sum_{0\leq j< N} \, \varphi(j)\cdot j.
\end{array}
\end{equation}

\noindent It is not hard to see that if $\|\varphi\| = K$, then
$\som(\varphi) \in \big\{0,\ldots, (N\minnetje 1)\cdot K\big\}$.

Accordingly, for $0 \leq i \leq (N\minnetje 1)\cdot K$ we define the
subset $\natMlt[K,i](\finset{N}) \subseteq \natMlt[K](\finset{N})$ as:
\[ \begin{array}{rcl}
\natMlt[K,i](\finset{N})
& \coloneqq &
\bigsetin{\varphi}{\natMlt[K](\finset{N})}{\som(\varphi) = i}.
\end{array} \]
\end{definition}

\auxproof{
See~\cite[Figure~5.3]{BenjaminQ03}: Distributing 10 kandies to 4
nerds $\{1,2,3,4\}$ can be done via a multiset $3\ket{1} + 2\ket{2}
+ 5\ket{4}$, as represented by the table:
\begin{center}
\begin{tabular}{ccccccccccccc}
$\star$ &
  $\star$ &
  $\star$ &
  $\mid$ &
  $\star$ &
  $\star$ &
  $\mid$ &
  $\mid$ &
  $\star$ &
  $\star$ &
  $\star$ &
  $\star$ &
  $\star$ 
\\
\hline
1 & 2 & 3 & 4 & 5 & 6 & 7 & 8 & 9 & 10 & 11 & 12 & 13
\end{tabular}
\end{center}
}

We also need the property of functoriality, which can be understood in
terms of repainting. Suppose we have a function $f\colon X \rightarrow
Y$, forming a transformation of $X$-colors to $Y$-colors. Then we can
turn it into a `repainting' function that takes an urn with
$X$-coloured balls into an urn with $Y$-coloured balls, by repainting
balls individually. The resulting function is written as $\natMlt(f)
\colon \natMlt(X) \rightarrow \natMlt(Y)$, where we thus apply
$\natMlt(-)$ both to sets and to functions between them. Usually this
does not lead to confusion. Explicitly,
\begin{equation}
\label{MltFunEqn}
\begin{array}{rcl}
\natMlt(f)\Big(\sum_{i} n_{i}\bigket{x_i}\Big)
& \coloneqq &
\sum_{i} n_{i}\bigket{f(x_{i})}.
\end{array} 
\end{equation}

\noindent It is not hard to see that identity and compositions are
preserved: $\natMlt(\idmap) = \idmap$ and $\natMlt(g \after f) =
\natMlt(g) \after \natMlt(f)$, and also that repainting maintains the
size of urns: $\big\|\natMlt(f)(\varphi)\big\| = \big\|\varphi\big\|$.

There is an obvious operation that turns sequences (tuples, lists)
into multisets by forgetting the order. We call it
\textit{accumulation} and write it as $\acc$. The above
urn~\eqref{Urn} can be obtained from a sequence via, for instance,
\[ \begin{array}{rcccl}
\acc\big(R,G,R,G,R,B\big)
& = &
3\bigket{R} + 2\bigket{G} + 1\bigket{B}
& = &
\acc\big(B,G,G,R,R,R\big).
\end{array} \]

\noindent For an arbitrary set $X$ and a number $K$ we can write
accumulation as a function $\acc \colon X^{K} \rightarrow
\natMlt[K](X)$, defined as $\acc\big(x_{1}, \ldots, x_{K}\big) \coloneqq
1\ket{x_1} + \cdots + 1\ket{x_K}$.

We include the following basic facts. Proofs may be found
in~\cite{Jacobs21g}.

\begin{lemma}
\label{MltLem}
Let $X$ be a finite set with $N \coloneqq \setsize{X} \geq 1$ elements.
\begin{enumerate}
\item \label{MltLemSize} There are $\big(\binom{N}{K}\big)$ multisets
  of size $K$ over $X$, that is, $\Bigsetsize{\natMlt[K](X)} =
  \displaystyle\bibinom{N}{K}$.

\item \label{MltLemAcc} For an arbitrary multiset
  $\varphi\in\natMlt[K](X)$, the number of sequences $\vec{x}\in X^{K}$
  with $\acc(\vec{x}) = \varphi$ is equal to the multiset coefficient
  $\coefm{\varphi}$, where:
\begin{equation}
\label{MltCoefEqn}
\begin{array}{rcl}
\coefm{\varphi}
& \coloneqq &
\displaystyle\frac{\|\varphi\|!}{\prod_{x} \varphi(x)!}.
\end{array} 
\end{equation}

\item \label{MltLemSum} For the sum of these multiset coefficients one
  has:
\begin{equation}
\label{MltSumEqn}
\begin{array}{rcl}
\displaystyle\sum_{\varphi\in\natMlt[K](X)} \, \coefm{\varphi}
& = &
N^{K}.
\end{array} 
\end{equation}
\end{enumerate}
\end{lemma}

Multinomial coefficients are often written of the form:
\[ \begin{array}{rcl}
\displaystyle\binom{n}{n_{1}, \ldots, n_{L}}
& = &
\displaystyle\frac{n!}{\prod_{i} n_{i}!}
   \qquad \mbox{where } \; n = \textstyle \sum_{i} n_{i}.
\end{array} \]

\noindent In this notation the order of the elements $n_i$ is irrelevant.
Hence it makes sense to define such coefficients for multisets. Using
multisets one can also give a snappy formulation of the multinomial
theorem, namely as:
\begin{equation}
\label{MultinomialThmEqn}
\begin{array}{rcl}
\Big(x_{0} + \cdots + x_{N-1}\Big)^{K}
& = &
\displaystyle\sum_{\varphi\in\natMlt[K](\finset{N})} \, 
   \coefm{\varphi}\cdot\prod_{0\leq j<N} x_{j}^{\varphi(j)},
\end{array}
\end{equation}

\noindent where, recall, $\finset{N} = \{0,\ldots,N\minnetje
1\}$.

In the sequel we shall see several forms of reversal, so we first
make this operation explicit.

\begin{lemma}
\label{RevertLem}
Consider the reversal self-isomorphism $\revert_{L} \colon
\{0,\ldots,L\} \congrightarrow \{0,\ldots,L\}$ by $\revert(i) =
L-i$. For a multiset $\varphi\in\natMlt[K](\finset{N})$ we use
functoriality to apply reversal to multisets:
\[ \begin{array}{rcl}
\Revert_{N-1}(\varphi)
\hspace*{\arraycolsep}\coloneqq\hspace*{\arraycolsep}
\natMlt(\revert_{N-1})(\varphi)
& = &
\displaystyle\sum_{0\leq i < N} \, \varphi(i)\bigket{\revert_{N-1}(i)}
\\[+1em]
& = &
\displaystyle\sum_{0\leq i < N} \, \varphi(i)\bigket{N\minnetje 1 - i}.
\end{array} \]

\noindent The following properties then hold.
\begin{enumerate}
\item The multiset coefficient of a reversed multiset is the same as
  the coefficient of the original: $\coefm{\Revert_{N-1}(\varphi)} =
  \coefm{\varphi}$. 

\item Using the sum of multisets $\som$ from Definition~\ref{SumDef},
  there is a commuting diagram:
\begin{equation}
\label{RevertSumDiag}
\vcenter{\xymatrix@R-0.8pc{
\natMlt[K](\finset{N})\ar[rr]^-{\som}\ar[d]_-{\Revert_{N-1}}^-{\cong} 
   & &
   \{0,1,\ldots,(N\minnetje 1)\cdot K\}\ar[d]^-{\revert_{(N-1)\cdot K}}_-{\cong}
\\
\natMlt[K](\finset{N})\ar[rr]^-{\som} & &
   \{0,1,\ldots,(N\minnetje 1)\cdot K\}
}}
\end{equation}

\noindent In particular, this gives: $\som(\varphi) = i \Longleftrightarrow
\som\big(\Revert_{K}(\varphi)\big) = (N\minnetje 1)\cdot K - i$.
\end{enumerate}
\end{lemma}

\begin{myproof}
\begin{enumerate}
\item This holds because the multiset coefficient is defined in
  terms of the multiplicities only --- see~\eqref{MltCoefEqn} ---
  which remain unchanged.

\item For a multiset $\varphi\in\natMlt[K](\finset{N})$ of size $K$,
\[ \begin{array}[b]{rcl}
\som\Big(\Revert_{N-1}(\varphi)\Big)
& = &
\displaystyle\sum_{0 \leq i < N} \, \varphi\big(N\minnetje 1-i\big)\cdot i
\\[+1.2em]
& = &
\displaystyle\sum_{0 \leq i < N} \, \varphi(i)\cdot\big(N\minnetje 1-i\big)
\\[+1.2em]
& = &
\displaystyle\sum_{0 \leq i < N} \, \varphi(i)\cdot (N\minnetje 1)
   - \sum_{0 \leq i < N} \, \varphi(i)\cdot i
\\[+0.6em]
& = &
K\cdot (N\minnetje 1) - \som(\varphi)
\hspace*{\arraycolsep}=\hspace*{\arraycolsep}
\revert_{(N-1)\cdot K}\Big(\som(\varphi)\Big).
\end{array} \eqno{\square} \]
\end{enumerate}
\end{myproof}

We now give an alternative formulation of the $N$-nomial
numbers~\eqref{NnomialSeqEqn}, not in terms of sequences (microstates)
but in terms of multisets (macrostates).

\begin{proposition}
\label{NnomialMltProp}
Let $N\geq 1$, $K\geq 0$ and $0 \leq i \leq (N\minnetje 1)\cdot K$ 
be given. 
\begin{enumerate}
\item \label{NnomialMltPropCoef} In terms of multisets one has:
\begin{equation}
\label{NnomialMltEqn}
\begin{array}{rcl}
C_{N}(K, i)
& = &
\displaystyle\sum_{\varphi\in\natMlt[K,i](\finset{N})} \, \coefm{\varphi}.
\end{array}
\end{equation}

\item \label{NnomialMltPropRev} $N$-nomials are closed under reversal:
\[ \begin{array}{rcl}
C_{N}(K,i)
& = &
C_{N}\Big(K, \, (N\minnetje 1)\cdot K - i\Big).
\end{array} \]

\noindent This generalises $\binom{K}{i} = \binom{K}{K-i}$ for $N=2$.
\end{enumerate}
\end{proposition}

\begin{myproof}
\begin{enumerate}
\item We use the accumulation function $\acc \colon \finset{N}^{K}
\rightarrow \natMlt[K](\finset{N})$, satisfying $\som(\vec{n}) = 
\som\big(\acc(\vec{n})\big)$ and $\bigsetsize{\acc^{-1}(\varphi)} = 
\coefm{\varphi}$, see Lemma~\ref{MltLem}~\eqref{MltLemAcc}.
\[ \begin{array}{rcl}
C_{N}(K, i)
\hspace*{\arraycolsep}\smash{\stackrel{\eqref{NnomialSeqEqn}}{=}}\hspace*{\arraycolsep}
\displaystyle\sum_{\vec{n} \in \finset{N}^{K}, \, \som(\vec{n}) = i} \, 1
& = &
\displaystyle\sum_{\varphi\in\natMlt[K](\finset{N})} \,
   \sum_{\vec{n} \in \acc^{-1}(\varphi), \, \som(\vec{n}) = i} \, 1
\\[+1.2em]
& = &
\displaystyle\sum_{\varphi\in\natMlt[K,i](\finset{N})} \,
   \sum_{\vec{n} \in \acc^{-1}(\varphi)} \, 1
\\[+1.2em]
& = &
\displaystyle\sum_{\varphi\in\natMlt[K,i](\finset{N})} \, 
  \coefm{\varphi}.
\end{array} \]

\item Now we can prove the closure of $N$-nomials under
reversal:
\[ \begin{array}[b]{rcl}
C_{N}(K, i)
& \smash{\stackrel{\eqref{NnomialMltEqn}}{=}} &
\displaystyle\sum_{\varphi\in\natMlt[K](\finset{N}), \, \som(\varphi) = i} \,
   \coefm{\varphi}
\\[+1.2em]
& = &
\displaystyle\sum_{\varphi\in\natMlt[K](\finset{N}), \, 
   \som(\natMlt(\revert_{N-1})(\varphi)) = i} \,
   \coefm{\natMlt(\revert_{N-1})(\varphi)}
\\[+1.2em]
& = &
\displaystyle\sum_{\varphi\in\natMlt[K](\finset{N}), \, 
   \som(\varphi) = (N-1)\cdot K - i} \, \coefm{\varphi}
\\[+1.2em]
& \smash{\stackrel{\eqref{NnomialMltEqn}}{=}} &
C_{N}\Big(K, \, (N\minnetje 1)\cdot K - i\Big).
\end{array} \eqno{\square} \]
\end{enumerate}
\end{myproof}

By combining the multiset formulation~\eqref{NnomialMltEqn} of
$N$-nomials with the multiset formulation of the Multinomial
Theorem~\eqref{MultinomialThmEqn} we can demonstrate that our approach
matches the definition of $N$-nomials in terms polynomial expressions,
occurring on the OEIS website~\cite{OEISNnomials}.

\begin{theorem}
\label{PolynomialThm}
For $N\geq 1$ and $K\geq 0$ one has, for an arbitrary variable $x$,
\[ \begin{array}{rcl}
\displaystyle\left(\sum_{0\leq j < N} \, x^{j}\right)^{K}
& = &
\displaystyle\sum_{0\leq i\leq (N-1)\cdot K} \, C_{N}(K,i)\cdot x^{i}.
\end{array} \]
\end{theorem}

\begin{myproof}
We use the Multinomial Theorem~\eqref{MultinomialThmEqn} in terms of
multisets in the first step in:
\[ \begin{array}[b]{rcl}
\displaystyle\left(\sum_{0\leq j < N} \, x^{j}\right)^{K}
& = &
\displaystyle\sum_{\varphi\in\natMlt[K](\finset{N})} \, 
   \coefm{\varphi}\cdot\prod_{0\leq j<N} \big(x^{j}\big)^{\varphi(j)}
\\[+1em]
& = &
\displaystyle\sum_{\varphi\in\natMlt[K](\finset{N})} \, 
   \coefm{\varphi}\cdot x^{\sum_{j} \varphi(j)\cdot j}
\\[+1em]
& = &
\displaystyle\sum_{\varphi\in\natMlt[K](\finset{N})} \, 
   \coefm{\varphi}\cdot x^{\som(\varphi)}
\\[+1em]
& = &
\displaystyle\sum_{0\leq i\leq (N-1)\cdot K} \, 
  \sum_{\varphi\in\natMlt[K](\finset{N}), \, \som(\varphi) = i} \, 
  \coefm{\varphi} \cdot x^{i}
\\[+1em]
& \smash{\stackrel{\eqref{NnomialMltEqn}}{=}} &
\displaystyle\sum_{0\leq i\leq (N-1)\cdot K} \, C_{N}(K,i)\cdot x^{i}.
\end{array} \eqno{\square} \]
\end{myproof}

\begin{remark}
\label{NnomialPythonRem}
We have given `mathematical' descriptions~\eqref{NnomialSeqEqn}
and~\eqref{NnomialMltEqn} of the $N$-nomial coefficient $C_{N}(K,i)$.
The sum of multiset coefficients in~\eqref{NnomialMltEqn} is a bit
more efficient to compute than the sequence version
in~\eqref{NnomialSeqEqn}, since it sums over all multisets, instead of
all sequences that accumulate to these multisets.  We like to add a
more efficient `computer science' description, in the programming
language \Python. It works by iteratively trying to collect the number
of sequences that add up to $i$.
\[ \vcenter{\hbox{\begin{minipage}[t]{0.6\textwidth}
\begin{lstlisting}[language=python]
def nomial(N, K, i):
    def collect(size, level):
        if size == 0:
            return 1 if level == 0 else 0
        if size == 1:
            return 1 if level < N else None
        out = 0
        for j in range(min(level+1, N)):
            col = collect(size - 1, level - j)
            if col != None: out += col
        return out if out > 0 else None
    return collect(K, i)
\end{lstlisting}
\end{minipage}}} \]
\end{remark}

\section{The Boltzmann distribution on multisets}\label{BoltzmannMltSec}

We have used ket notation $\ket{-}$ to separate multiplicities and
elements in multisets. It makes sense to use the same notation for
(finite, discrete) probability distributions, now to separate
probabilities and elements. The distribution corresponding to the urn
$\upsilon$ in~\eqref{Urn} is then written as $\frac{1}{2}\ket{R} +
\frac{1}{3}\ket{G} + \frac{1}{6}\ket{B}$. It gives the probabilities
of randomlingly drawing a ball from $\upsilon$, per colour. We briefly
discuss (discrete) distributions in general, before concentrating
on Boltzmann distributions.

We shall write $\Dst(X)$ for the set of distributions over a set $X$.
The elements of $\Dst(X)$ may be written in ket form as finite formal
convex sums $\sum_{i} r_{i}\ket{x_i}$, where $x_{i} \in X$ and
$r_{i}\in [0,1]$ with $\sum_{i} r_{i} = 1$. Equivalently, such
distributions can be written in functional form as functions
$\omega\colon X \rightarrow [0,1]$ with finite support
$\support(\omega) \coloneqq \setin{x}{X}{\omega(x) \neq 0}$ and with
$\sum_{x} \omega(x) = 1$.

There is a systematic way to turn a (non-empty) multiset into
a distribution, namely by normalisation. We write this operation
as $\flrn$ and call it frequentist learning, since it involves learning
a distribution by counting. Explicitly, when $\varphi \neq \zero$,
that is, when $\|\varphi\| \neq 0$,
\begin{equation}
\label{FlrnEqn}
\begin{array}{rcl}
\flrn(\varphi)
& \coloneqq &
\displaystyle\sum_{x\in\support(\varphi)} \, \frac{\varphi(x)}{\|\varphi\|}\,
   \bigket{x}.
\end{array}
\end{equation}

\noindent For instance, $\flrn\big(3\ket{a} + 4\ket{b} + 5\ket{c}\big)
= \frac{1}{4}\ket{a} + \frac{1}{3}\ket{b} + \frac{5}{12}\ket{c}$.

\begin{definition}
\label{PushDef}
Let $\omega = \sum_{i} r_{i}\ket{x_i} \in\Dst(X)$ be a distribution on
a set $X$, together with another set $Y$.
\begin{enumerate}
\item For a function $f\colon X \rightarrow Y$ there is the
\emph{image} distribution $\Dst(f)(\omega) \in \Dst(Y)$,
defined as in~\eqref{MltFunEqn}:
\begin{equation}
\label{DstFunEqn}
\begin{array}{rcl}
\Dst(f)\Big(\sum_{i} r_{i}\bigket{x_}\Big) 
& \coloneqq &
\sum_{i} r_{i}\bigket{f(x_{i})}.
\end{array}
\end{equation}

\item For a `channel' $c\colon X \rightarrow \Dst(Y)$ there is the
  \emph{pushforward} distribution $c_{*}(\omega) \in \Dst(Y)$ given
  by:
\begin{equation}
\label{PushEqn}
\begin{array}{rcccl}
c_{*}\Big(\sum_{i} r_{i}\bigket{x_}\Big) 
& \coloneqq &
\sum_{i} r_{i}\cdot c(x_{i})
& = &
\displaystyle\sum_{y\in Y} \, \textstyle 
   \Big(\sum_{i} r_{i}\cdot c(x_{i})(y)\Big)\,\bigket{y}.
\end{array}
\end{equation}

\noindent For another channel $d\colon Y \rightarrow \Dst(Z)$ we can
now define the channel composite $d \klafter c \colon X \rightarrow
\Dst(Z)$ as $\big(d \klafter c\big)(x) = d_{*}\big(c(x)\big)$.
\end{enumerate}
\end{definition}

A channel $c\colon X \rightarrow \Dst(Y)$ gives a distribution $c(x)
\in \Dst(Y)$ for each $x\in X$. This corresponds to what is commonly
described as a conditional probability $\Prob(y|x)$.  The pushforward
can be understood as application of the channel `in probability', as
expressed by the formulation in the middle~\eqref{PushEqn}.

We shall use distributions over numbers, but also distributions over
multisets. The latter will be written with kets-inside-kets. For
instance, the equation $\sum_{\varphi\in\natMlt[K](X)} \coefm{\varphi}
= N^{K}$ in Lemma~\ref{MltLem}~\eqref{MltLemSum} gives rise to
a \emph{multiset coefficient} distribution of the form:
\begin{equation}
\label{MltCoefDstEqn}
\begin{array}{rcl}
\multisetcoefficientdistribution[X,K]
& = &
\displaystyle\sum_{\varphi\in\natMlt[K](X)} \, \frac{\coefm{\varphi}}{N^{K}}
   \,\bigket{\varphi}
\qquad\mbox{where } N = \setsize{X}.
\end{array}
\end{equation}

\noindent This distribution is identified in~\cite{JacobsS25} as a
stationary distribution for a certain Markov chain on $\natMlt[K](X)$.
We give an illustration of this distribution~\eqref{MltCoefDstEqn},
for the set $X = \finset{3} = \{0,1,2\}$ and for size $K=2$. We then
get a distribution over the $\big(\binom{3}{2}\big) = 6$ multisets in
$\natMlt[2]\big(\finset{3}\big)$ of the form:
\[ \begin{array}{rcl}
\multisetcoefficientdistribution[\finset{3},2]
& = &
\frac{1}{9}\Bigket{2\ket{0}} + 
\frac{2}{9}\Bigket{1\ket{0} + 1\ket{1}} + 
\frac{1}{9}\Bigket{2\ket{1}}
\\[+0.2em]
& & \qquad +\;
\frac{2}{9}\Bigket{1\ket{0} + 1\ket{2}} + 
\frac{2}{9}\Bigket{1\ket{1} + 1\ket{2}} + 
\frac{1}{9}\Bigket{2\ket{2}}.
\end{array} \]



\noindent We write this distribution via nested kets, with the outer
big kets for the distribution (with probabilities in front), and
inside the big kets the multisets of size $2$ over the set $\finset{3}
= \{0,1,2\}$, expressed via the smaller kets.

We are not interested in the entire multiset coefficient distribution
in~\eqref{MltCoefDstEqn} but in a restriction to those multisets with
a fixed sum $i$, as in Definition~\ref{SumDef}. In order to turn these
into a proper distribution --- with probabilities adding up to one ---
we have to normalise, precisely with the $N$-nomial coefficients,
using Proposition~\ref{NnomialMltProp}~\eqref{NnomialMltPropCoef}.
This leads to our first Boltzmann distribution.

\begin{definition}
\label{BoltzmannMltDef}
Fix numbers $N\geq 1$ and $K\geq 0$. The \emph{Boltzmann-on-multisets}
distributions arise via the function $\boltzmannmlt[N,K] \colon
\big\{0,\ldots,(N\minnetje 1)\cdot K\big\} \rightarrow
\Dst\Big(\natMlt[K](\finset{N})\Big)$ defined as:
\begin{equation}
\label{BoltzmannMltEqn}
\begin{array}{rcl}
\boltzmannmlt[N, K](i)
& \coloneqq &
\displaystyle\sum_{\varphi\in\natMlt[K,i](\finset{N})} \, 
   \frac{\coefm{\varphi}}{C_{N}(K,i)} \, \bigket{\varphi}.
\end{array}
\end{equation}
\end{definition}

As an aside, we put the parameters, like $N$ and $K$
in~\eqref{BoltzmannMltEqn}, between square brackets $[-]$ when they
determine the type of the operation at hand. The parameter(s) between
round brackets $(-)$ are inputs to a
function\textsuperscript{1}\footnotetext{\textsuperscript{1}
  One could narrow down the type of the Boltzmann-on-multisets
  distribution to a dependent function of the form $\boltzmannmlt[N,K]
  \colon \prod_{0\leq i\leq (N-1)\cdot K}
  \Dst\big(\natMlt[K,i](\finset{N})\big)$.}.

We illustrate how this Boltzmann-on-multisets distribution arises in
physics. The set $\finset{N} = \{0,\ldots,N-1\}$ captures the
different energy levels at which particles can be, where $K$ is the
number of particles.  A multiset $\varphi\in\natMlt[K,i](\finset{N})$
with $\som(\varphi) = i$ is then a configuration (`composition'
in~\cite{DillB10} or `macrostate' in~\cite{PathriaB11}) of $K$
particles, at various energy levels, so that the total energy
$\som(\varphi) = \sum_{j} \varphi(j)\cdot j$ is equal to $i$, for all
configurations under consideration.

\begin{example}
\label{PhysicsMltEx}
In~\cite[Appendix~C]{EisbergR85} an example is given with $N=4$ energy
levels, with $K=4$ particles, and with total energy $i=3$.  There are
three multisets in $\varphi\in\natMlt[3](\finset{4})$ with
$\som(\varphi) = 3$. The following table uses the particle
configurations from~\cite{EisbergR85} in the column on the left, and
the fractions in the column on right. Inbetween, these configurations
and fractions are interpreted in the current setting with multisets
and their coefficients. The checkmarks $\checkmark$ indicate how many
particles are at which energy level.
\begin{center}
\begin{tabular}{cccc}
\textbf{configuration over} $0,\ldots, 4$ & 
   \textbf{multiset} $\varphi$ & 
   \hspace*{-0.5em}\textbf{coefficient} $\coefm{\varphi}$\hspace*{-0.5em} & 
   \textbf{fraction}
\\
\hline
\hline
\\[-1.0em]
$\vcenter{\hbox{\includegraphics[scale=0.11]{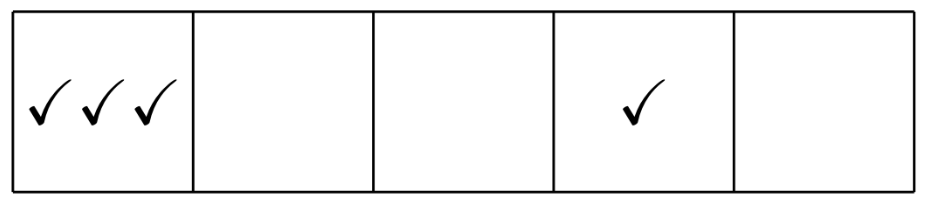}}}$
   & $3\ket{0} + 1\ket{3}$ & 4 & 
   $\displaystyle\frac{4}{4 \plusje 12 \plusje 4} = \frac{1}{5}$
\\[+1.5em]
$\vcenter{\hbox{\includegraphics[scale=0.11]{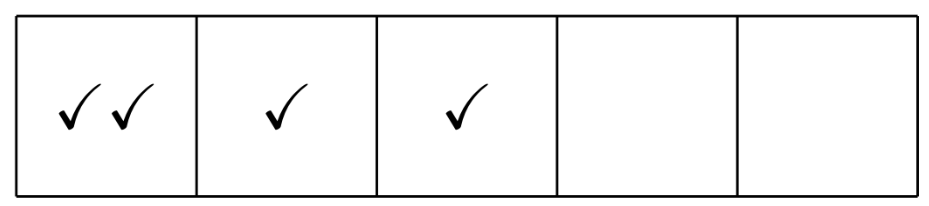}}}$
   & $2\ket{0} + 1\ket{1} + 1\ket{2}$ & 12 & 
   $\displaystyle\frac{12}{4 \plusje 12 \plusje 4} = \frac{3}{5}$
\\[+1.5em]
$\vcenter{\hbox{\includegraphics[scale=0.11]{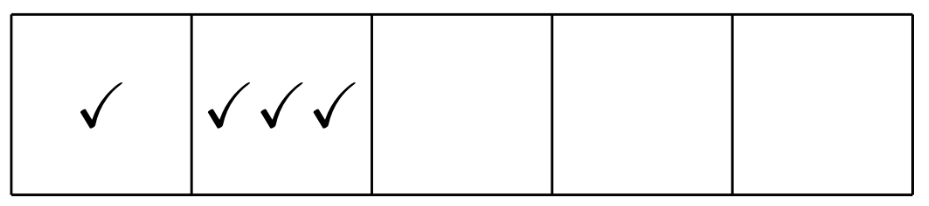}}}$
   & $1\ket{0} + 3\ket{1}$ & 4 & 
   $\displaystyle\frac{4}{4 \plusje 12 \plusje 4} = \frac{1}{5}$
\end{tabular}
\end{center}

\noindent The resulting Boltzmann-on-multisets distribution indeed
contains the fractions in the last column as probabilities, since the
$4$-nomial coefficient $C_{4}(4,3)$ is equal to~$20$. This is by
construction, the sum of the coefficients in the above
coefficient-column. Then:
\begin{equation}
\label{BoltzmannMltEisberg}
\begin{array}{rcl}
\lefteqn{\boltzmannmlt[4, 4](3)
\hspace*{\arraycolsep}=\hspace*{\arraycolsep}
\displaystyle\sum_{\varphi\in\natMlt[4, 3](\finset{4})} \,
    \frac{\coefm{\varphi}}{C_{4}(4,3)}\,\bigket{\varphi}}
\\[+1.2em]
& = &
\frac{1}{5}\Bigket{3\ket{0} + 1\ket{3}} + 
   \frac{3}{5}\Bigket{2\ket{0} + 1\ket{1} + 1\ket{2}} +
   \frac{1}{5}\Bigket{1\ket{0} + 3\ket{1}}.
\end{array}
\end{equation}

\noindent The book~\cite{EisbergR85} describes this situation as an
illustration, but does not contain the general
distribution~\eqref{BoltzmannMltEqn}, nor the $N$-nomial numbers that
are needed for normalisation.  In fact, it does not mention multisets
at all.

A similar but bigger table, like the one above, can be found in the
online resources~\cite{TiplerLResources} for students and teachers of
the book~\cite{TiplerL12}. It builds up the Boltzmann-on-multisets
distribution $\boltzmannmlt[9,6](8)$, with $N=9$ energy levels
$0,\ldots,8$, with $K=6$ particles, having total energy $i=8$. This
leads to a distribution on $20$~multisets, that are listed there, with
sum of multinomial coefficients $C_{9}(6,8) = 1287$. Further details
about this illustration are provided in the appendix.

\ignore{

}

\end{example}

We conclude this section with two observations about the
Boltzmann-on-multisets distribution, namely that it is stable under
reversal and that it can be obtained as image distribution from a
uniform distribution on sequences (microstates). In order to express
this we need the functoriality of both $\Dst$ and $\natMlt$.

\begin{lemma}
\label{BoltzmannMltLem}
Fix numbers $N \geq 1$ and $K\geq 0$.
\begin{enumerate}
\item \label{BoltzmannMltLemRev} The Boltzmann-on-multisets
  distribution~\eqref{BoltzmannMltEqn} is stable under reversal as
  expressed by the following commuting diagram, using the revert
  isomorphism from Lemma~\ref{RevertLem}.
\[ \xymatrix@R-0.8pc@C+1pc{
\{0,1,\ldots,(N\minnetje 1)\cdot K\}\ar[d]_-{\revert_{(N-1)\cdot K}}^-{\cong}
   \ar[rr]^-{\boltzmannmlt[N, K]} & &
   \Dst\Big(\natMlt[K](\finset{N})\Big)
   \ar[d]^-{\Dst\natMlt(\revert_{N-1})}_-{\cong}
\\
\{0,1,\ldots,(N\minnetje 1)\cdot K\}\ar[rr]^-{\boltzmannmlt[N, K]} & &
   \Dst\Big(\natMlt[K](\finset{N})\Big)
} \]

\item \label{BoltzmannMltLemIm} For a number $0 \leq i \leq
  (N\minnetje 1)\cdot K$, consider the set of sequences $S$ defined
  below, as in~\eqref{NnomialSeqEqn}, with the uniform distribution
  $\unif$ on it.
\[ \begin{array}{rclcrcl}
S
& \coloneqq &
\bigsetin{v}{\finset{N}^{K}}{\som(v) = i}
& \qquad\qquad &
\unif
& \coloneqq &
\displaystyle\sum_{v\in S} \, \frac{1}{C_{N}(K,i)}\,\bigket{v}.
\end{array} \]

\noindent Then: $\Dst(\acc)(\unif) = \boltzmannmlt[N,K](i)$.
\end{enumerate}
\end{lemma}

\begin{myproof}
\begin{enumerate}
\item We use the closure of $N$-nomials under reversal from
Proposition~\ref{NnomialMltProp}~\eqref{NnomialMltPropRev}, and also
the commutation of $\revert$ with sums~\eqref{RevertSumDiag},
where we wrote $\Revert_{N-1}(\varphi) = \natMlt(\revert_{N-1})(\varphi)$.
\[ \begin{array}[b]{rcl}
\lefteqn{\Dst\natMlt(\revert_{N-1})\Big(\boltzmannmlt[N,K](i)\Big)}
\\[0.2em]
& \smash{\stackrel{\eqref{BoltzmannMltEqn}}{=}} &
\displaystyle\sum_{\varphi\in\natMlt[K](\finset{N}), \, \som(\varphi) = i} \,
   \frac{\coefm{\varphi}}{C_{N}(K,i)}\,\bigket{\natMlt(\revert_{N-1})(\varphi)}
\\[+1em]
& = &
\displaystyle\sum_{\varphi\in\natMlt[K](\finset{N}), \, 
   \som(\natMlt(\revert_{N-1})(\varphi)) = i} \,
   \frac{\coefm{\natMlt(\revert_{K})(\varphi)}}{C_{N}(K,i)}\,\bigket{\varphi}
\\[+1em]
& = &
\displaystyle\sum_{\varphi\in\natMlt[K](\finset{N}), \, \som(\varphi) = (N-1)\cdot K - i} \,
   \frac{\coefm{\varphi}}{C_{N}(K, (N\minnetje 1)\cdot K - i)}\,\bigket{\varphi}
\\[+1.2em]
& = &
\boltzmannmlt[N,K]\Big((N\minnetje 1)\cdot K - i\Big)
\hspace*{\arraycolsep}=\hspace*{\arraycolsep}
\boltzmannmlt[N,K]\Big(\revert_{(N\minnetje 1)\cdot K}(i)\Big).
\end{array} \]

\item Since:
\[ \begin{array}[b]{rcl}
\Dst(\acc)(\unif)
& = &
\displaystyle\sum_{v\in S} \, \frac{1}{C_{N}(K,i)} \, \bigket{\acc(v)}
\\
& = &
\displaystyle\sum_{\varphi\in\natMlt[K,i](\finset{N})} \;
   \sum_{v\in\acc^{-1}(\varphi)} \, \frac{1}{C_{N}(K,i)} \, \bigket{\acc(v)}
\\[+1em]
& = &
\displaystyle\sum_{\varphi\in\natMlt[K,i](\finset{N})} \;
   \sum_{v\in\acc^{-1}(\varphi)} \, \frac{1}{C_{N}(K,i)} \, \bigket{\varphi}
\\[+1em]
& = &
\displaystyle\sum_{\varphi\in\natMlt[K,i](\finset{N})} \,
   \frac{\coefm{\varphi}}{C_{N}(K,i)} \, \bigket{\varphi}
   \qquad \mbox{by Lemma~\ref{MltLem}~\eqref{MltLemAcc}}
\\[+0.8em]
& = &
\boltzmannmlt[N,K](i)
\end{array} \eqno{\square} \]

\end{enumerate}
\end{myproof}

\ignore{

}

\section{The Boltzmann distribution on numbers}\label{BoltzmannNumSec}

We now transform the Boltzmann-on-multisets distribution
$\boltzmannmlt$ into a `Boltzmann-on-numbers' distribution
$\boltzmannnum$. In physical terms, it describes the spread of
particles over the various energy levels $\{0,\ldots,N\minnetje
1\}$. There is a systematic way to do this transformation, via
frequentist learning $\flrn$ from~\eqref{FlrnEqn}. The triangle below
defines $\boltzmannnum$ as composition of two functions.
\begin{equation}
\label{BoltzmannNumDiag}
\vcenter{\xymatrix@R-1pc{
& & \Dst\Big(\natMlt[K](\finset{N})\Big)\ar[dd]^-{\flrn_{*}}
\\
\big\{0,\ldots,(N\minnetje 1)\cdot K\big\}
   \ar@/^3ex/[urr]^-{\boltzmannmlt[N,K]}
   \ar@/_3ex/[drr]_-{\boltzmannnum[N,K]} & &
\\
& & \Dst\big(\finset{N}\big)\rlap{$\;= \Dst\big(\{0,\ldots,N\minnetje 1\}\big)$}
}}\hspace*{5em}
\end{equation}

\noindent The diagram contains the pushforward $\flrn_{*} \colon
\Dst\big(\natMlt[K](X)\big) \rightarrow \Dst(X)$ of the frequentist
learning channel $\flrn \colon \natMlt[K](X) \rightarrow \Dst(X)$, as
introduced in~\eqref{PushEqn}. As a consequence,
$\boltzmannnum[N,K](i)$ gives the energy for a randomly chosen
particle.

The next definition gives both an abstract and a concrete formulation
of the Boltz\-mann-on-numbers distribution, so the above triangle serves
only as background information.

\begin{definition}
\label{BoltzmannNumDef}
For numbers $N,K\geq 1$ and $0\leq i\leq (N\minnetje 1)\cdot K$ we
define the Boltzmann-on-numbers distribution as:
\begin{equation}
\label{BoltzmannNumEqn}
\begin{array}{rcl}
\boltzmannnum[N, K](i)
& \coloneqq &
\flrn_{*}\Big(\boltzmannmlt[N, K](i)\Big)
\\[+0.6em]
& = &
\displaystyle\sum_{0\leq j < N} \, 
   \sum_{\varphi\in\natMlt[K,i](\finset{N})} \,
   \frac{\coefm{\varphi}\cdot \varphi(j)}{C_{N}(K,i) \cdot K}
   \,\bigket{j}.
\end{array}
\end{equation}

\noindent A subtle point is that if $i < N$, then we can limit the sum
over $j$ to $j \leq i$, since if $\som(\varphi) = i$, then $\varphi(j)
= 0$ for $j > i$. We can make this explicit by letting $j$ range over
$0 \leq j < \min(N, i+1)$. The next section concentrates on this
special case $i < N$.
\end{definition}

We double-check that the probabilities in~\eqref{BoltzmannNumEqn} add
up to one:
\[ \begin{array}{rcl}
\displaystyle \sum_{0\leq j < N} \, 
   \sum_{\varphi\in\natMlt[K,i](\finset{N})} \,
   \frac{\coefm{\varphi}\cdot \varphi(j)}{C_{N}(K,i) \cdot K}
& = &
\displaystyle 
   \sum_{\varphi\in\natMlt[K,i](\finset{N})} \,
   \frac{\coefm{\varphi}\cdot \sum_{0\leq j < N}\varphi(j)}{C_{N}(K,i) \cdot K}
\\[+1em]
& = &
\displaystyle \sum_{\varphi\in\natMlt[K,i](\finset{N})} \,
   \frac{\coefm{\varphi}\cdot \|\varphi\|}{C_{N}(K,i) \cdot K}
\\[+1em]
& = &
\displaystyle \sum_{\varphi\in\natMlt[K,i](\finset{N})} \,
   \frac{\coefm{\varphi}}{C_{N}(K,i)}
\hspace*{\arraycolsep}\smash{\stackrel{\eqref{NnomialMltEqn}}{=}}\hspace*{\arraycolsep}
1.
\end{array} \]

\begin{example}
\label{PhysicsNumEx}
In the context of Example~\ref{PhysicsMltEx}, the associated
Boltzmann-on-numbers distribution will be computed in two ways.
First we shall use the formulation via the extension
$\flrn_*$ of frequentist learning --- as used in the
middle of~\eqref{BoltzmannNumEqn} --- to illustrate how this works.
\[ \begin{array}{rcl}
\lefteqn{\boltzmannnum[4, 4](3)}
\\
& \smash{\stackrel{\eqref{BoltzmannNumEqn}}{=}} &
\flrn_{*}\Big(\boltzmannmlt[4, 4](3)\Big)
\\
& \smash{\stackrel{\eqref{BoltzmannMltEisberg}}{=}} &
\flrn_{*}\Big(\frac{1}{5}\Bigket{3\ket{0} + 1\ket{3}} + 
   \frac{3}{5}\Bigket{2\ket{0} + 1\ket{1} + 1\ket{2}} +
   \frac{1}{5}\Bigket{1\ket{0} + 3\ket{1}}\Big)
\\[+0.2em]
& \smash{\stackrel{\eqref{PushEqn}}{=}} &
\frac{1}{5}\cdot\flrn\Big(3\ket{0} + 1\ket{3}\Big) + 
   \frac{3}{5}\cdot \flrn\Big(2\ket{0} + 1\ket{1} + 1\ket{2}\Big) +
   \frac{1}{5}\cdot \flrn\Big(1\ket{0} + 3\ket{1}\Big)
\\[+0.2em]
& = &
\frac{1}{5}\cdot\Big(\frac{3}{4}\ket{0} + \frac{1}{4}\ket{3}\Big) + 
   \frac{3}{5}\cdot \Big(\frac{1}{2}\ket{0} + 
      \frac{1}{4}\ket{1} + \frac{1}{4}\ket{2}\Big) +
   \frac{1}{5}\cdot \Big(\frac{1}{4}\ket{0} + \frac{3}{4}\ket{1}\Big)
\\[+0.2em]
& = &
\Big(\frac{1}{5}\cdot\frac{3}{4} + 
   \frac{3}{5}\cdot \frac{1}{2} + 
   \frac{1}{5}\cdot \frac{1}{4}\Big)\ket{0} +
\Big(\frac{3}{5}\cdot \frac{1}{4} + 
   \frac{1}{5}\cdot \frac{3}{4}\Big)\ket{1} +
\frac{3}{5}\cdot \frac{1}{4}\ket{2} + 
\frac{1}{5}\cdot \frac{1}{4}\ket{3}
\\[+0.2em]
& = &
\frac{1}{2}\bigket{0} + \frac{3}{10}\bigket{1} + 
   \frac{3}{20}\bigket{2} + \frac{1}{20}\bigket{3}.
\end{array} \]

\noindent We also elaborate the formulation on the
right-hand-side of~\eqref{BoltzmannNumEqn}.
\[ \begin{array}{rcl}
\lefteqn{\boltzmannnum[4, 4](3)}
\\
& \smash{\stackrel{\eqref{BoltzmannNumEqn}}{=}} &
\displaystyle\sum_{0\leq j < 4} \, 
   \sum_{\varphi\in\natMlt[4,3](\finset{4})} \,
   \frac{\coefm{\varphi}\cdot \varphi(j)}{C_{N}(K,i) \cdot K}
   \,\bigket{j}
\\[+1.4em]
& = &
\displaystyle\frac{4\cdot 3 + 12\cdot 2 + 4\cdot 1}{4\cdot 20}\bigket{0}
   + \frac{12\cdot 1 + 4\cdot 3}{4\cdot 20}\bigket{1}
   + \frac{12\cdot 1}{4\cdot 20}\bigket{2}
   + \frac{4\cdot 1}{4\cdot 20}\bigket{3}
   + \frac{0}{4\cdot 20}\bigket{4}
\\[+0.6em]
& = &
\frac{1}{2}\bigket{0} + \frac{3}{10}\bigket{1} + 
   \frac{3}{20}\bigket{2} + \frac{1}{20}\bigket{3}.
\end{array} \]

\noindent The calculations in the second line, matching the table in
Example~\ref{PhysicsMltEx}, are performed concretely
in~\cite[Appendix~C]{EisbergR85}, but without a general formula, as
in~\eqref{BoltzmannNumEqn}.
\end{example}

\begin{figure}
\begin{center}
\begin{tabular}{cccc}
\includegraphics[scale=0.13]{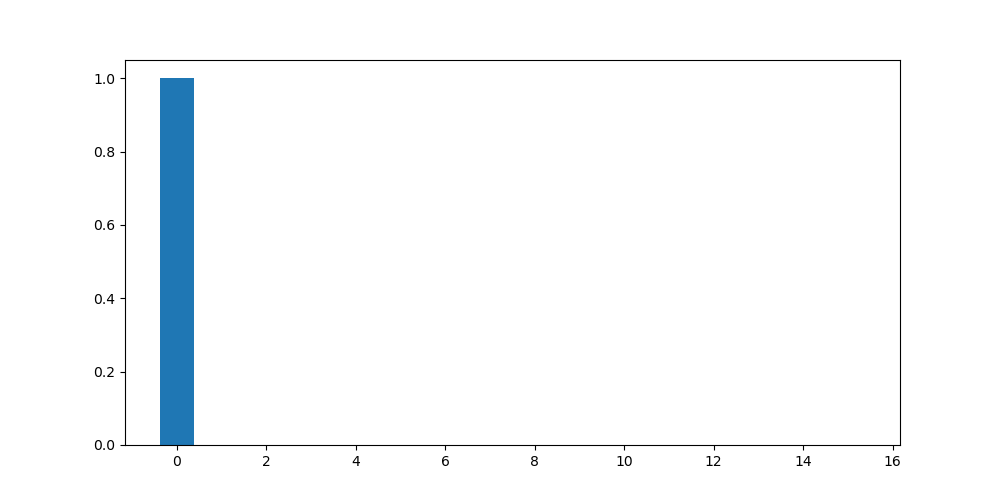}
   \hspace*{-1em} & \hspace*{-1em}
\includegraphics[scale=0.13]{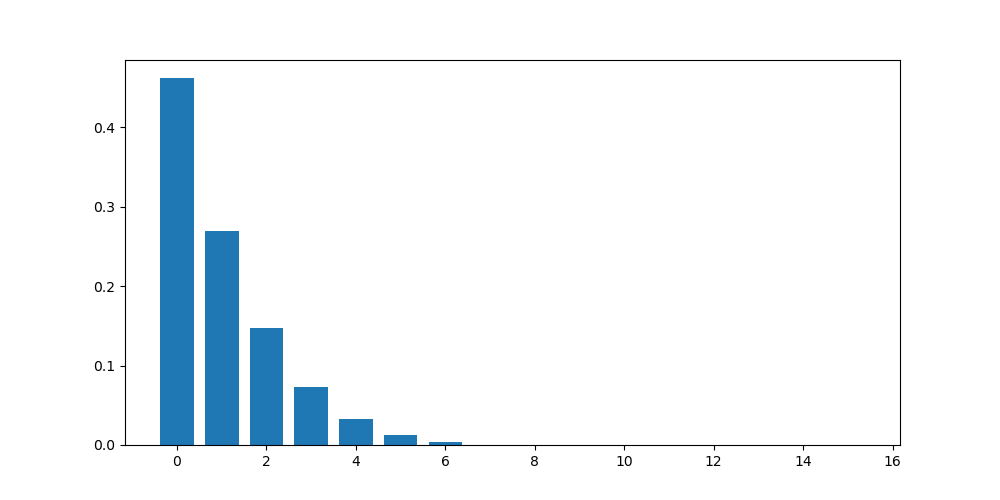}
   \hspace*{-1em}  & \hspace*{-1em}
\includegraphics[scale=0.13]{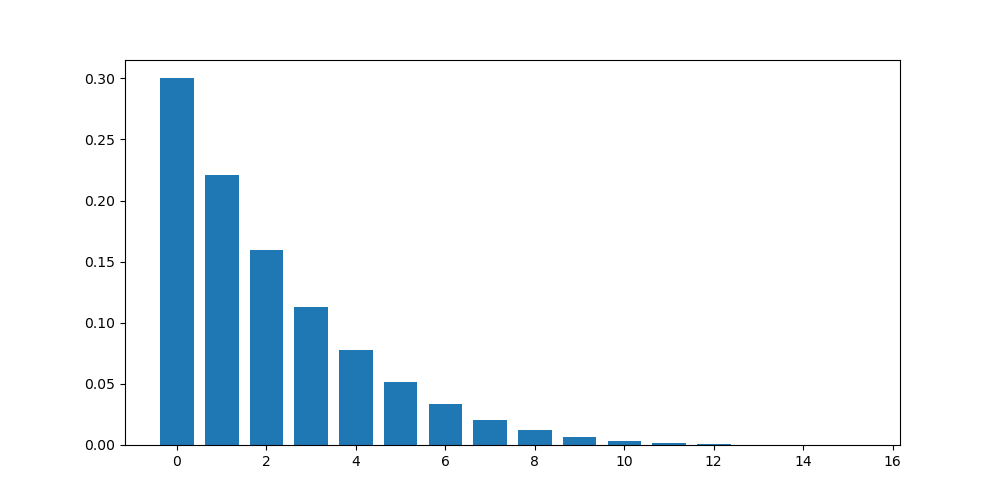}
   \hspace*{-1em}  & \hspace*{-1em}
\includegraphics[scale=0.13]{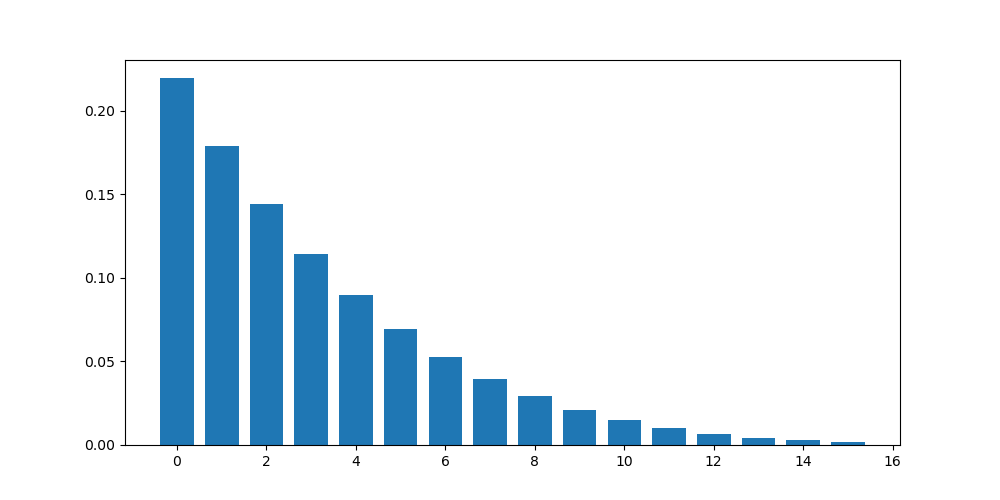}
\\
\includegraphics[scale=0.13]{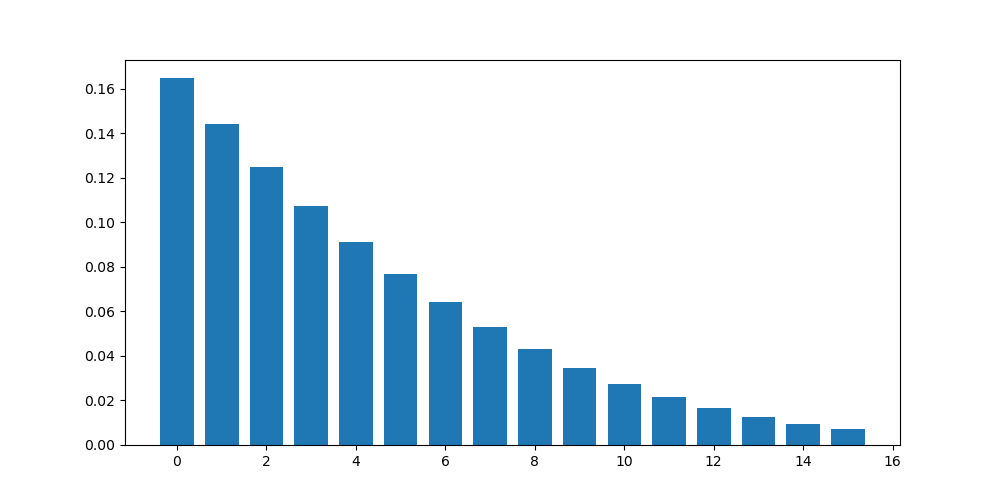}
   \hspace*{-1em}  & \hspace*{-1em}
\includegraphics[scale=0.13]{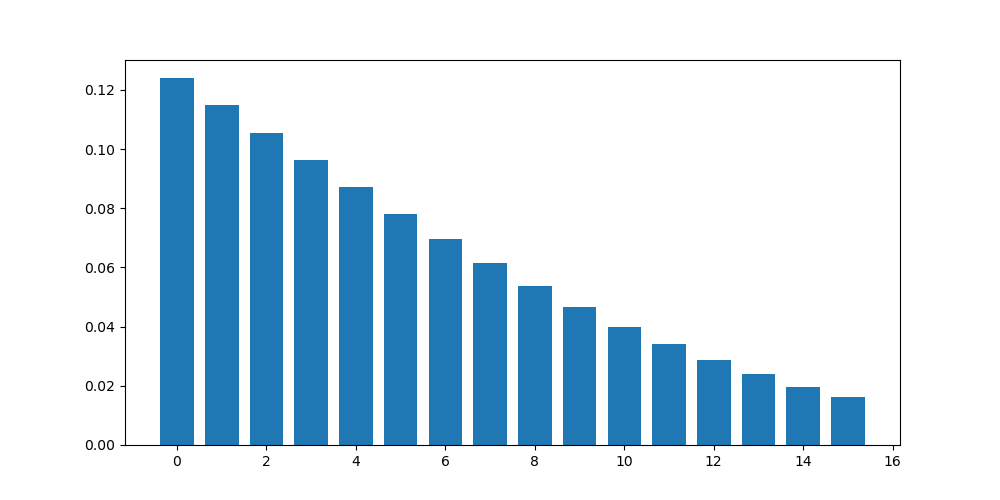}
   \hspace*{-1em}  & \hspace*{-1em}
\includegraphics[scale=0.13]{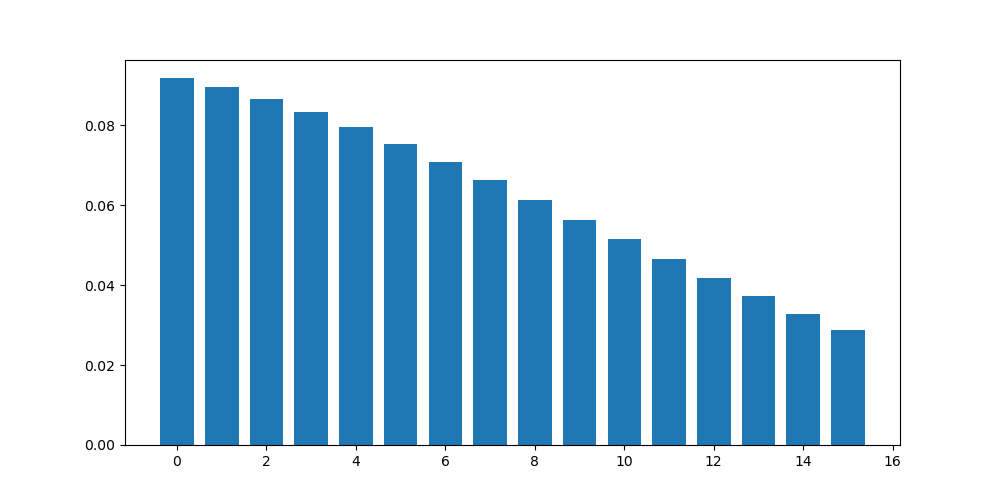}
   \hspace*{-1em}  & \hspace*{-1em}
\includegraphics[scale=0.13]{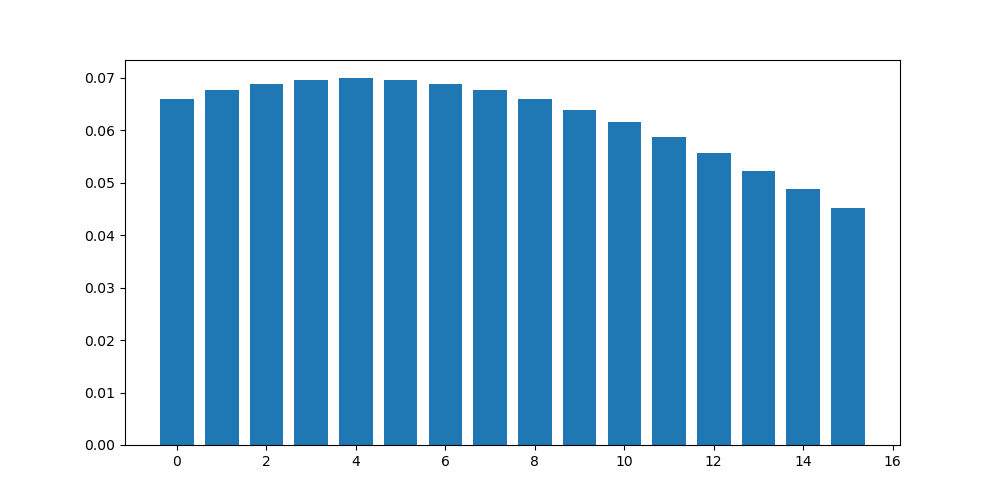}
\\
\includegraphics[scale=0.13]{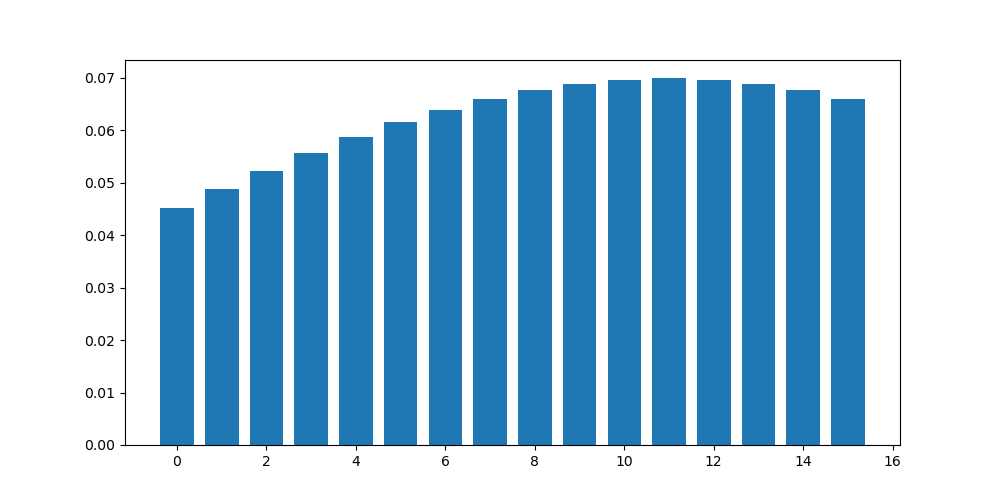}
   \hspace*{-1em}  & \hspace*{-1em}
\includegraphics[scale=0.13]{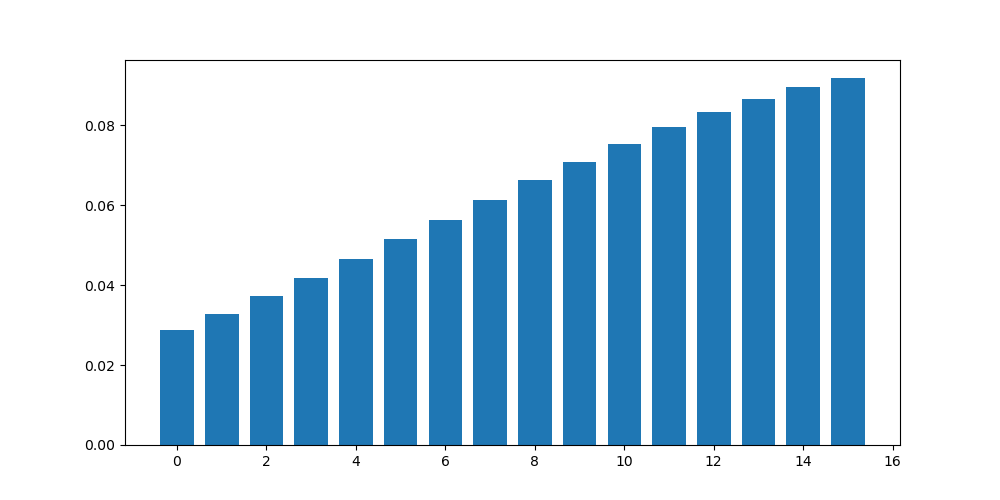}
   \hspace*{-1em}  & \hspace*{-1em}
\includegraphics[scale=0.13]{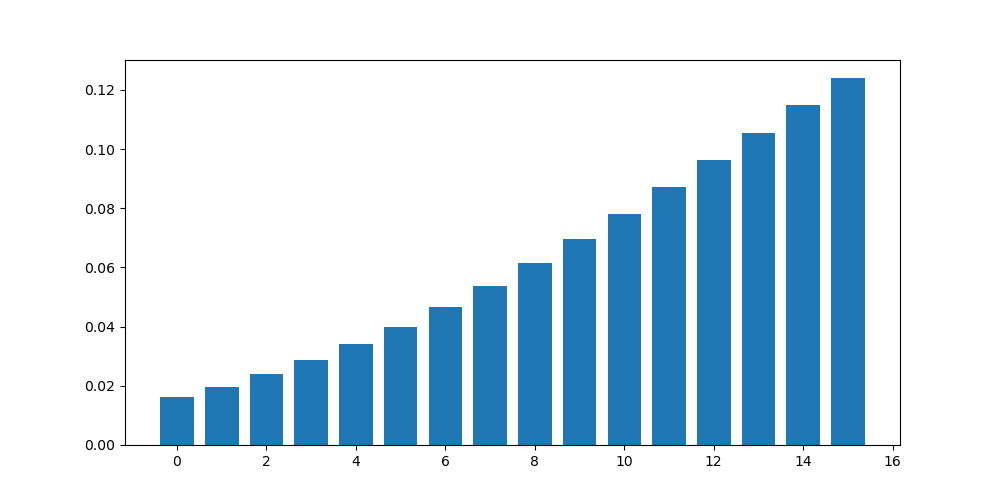}
   \hspace*{-1em}  & \hspace*{-1em}
\includegraphics[scale=0.13]{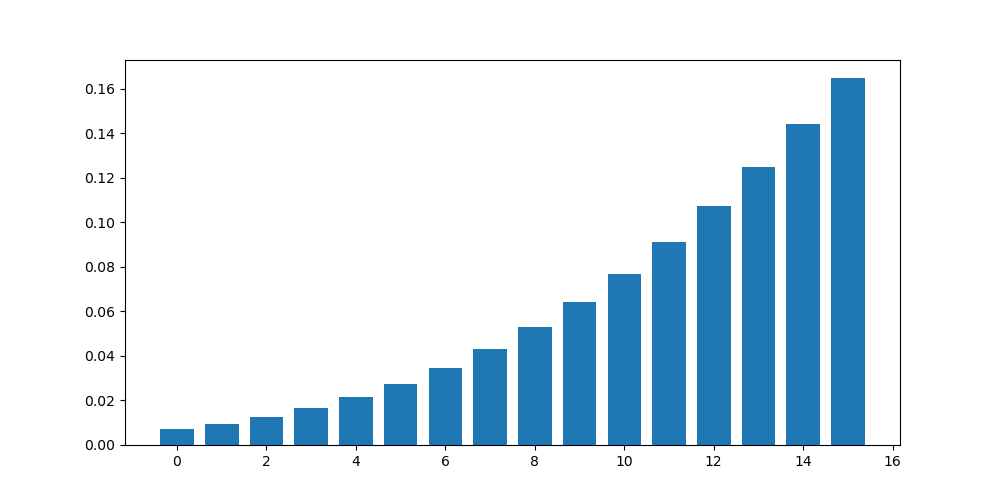}
\\
\includegraphics[scale=0.13]{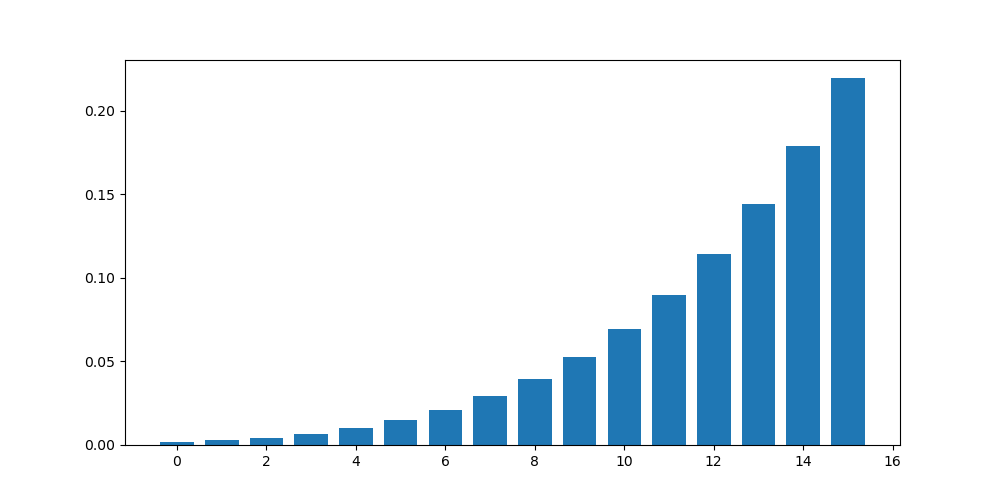}
   \hspace*{-1em}  & \hspace*{-1em}
\includegraphics[scale=0.13]{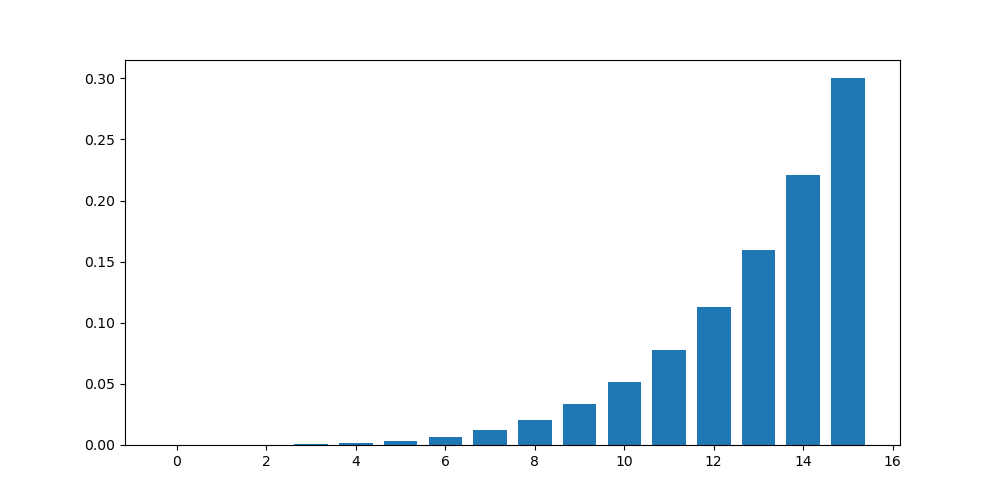}
   \hspace*{-1em}  & \hspace*{-1em}
\includegraphics[scale=0.13]{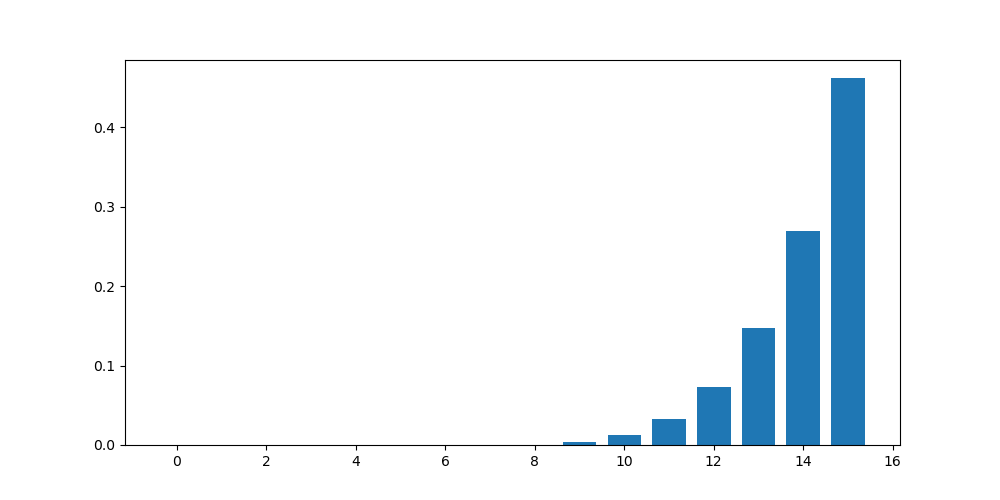}
   \hspace*{-1em}  & \hspace*{-1em}
\includegraphics[scale=0.13]{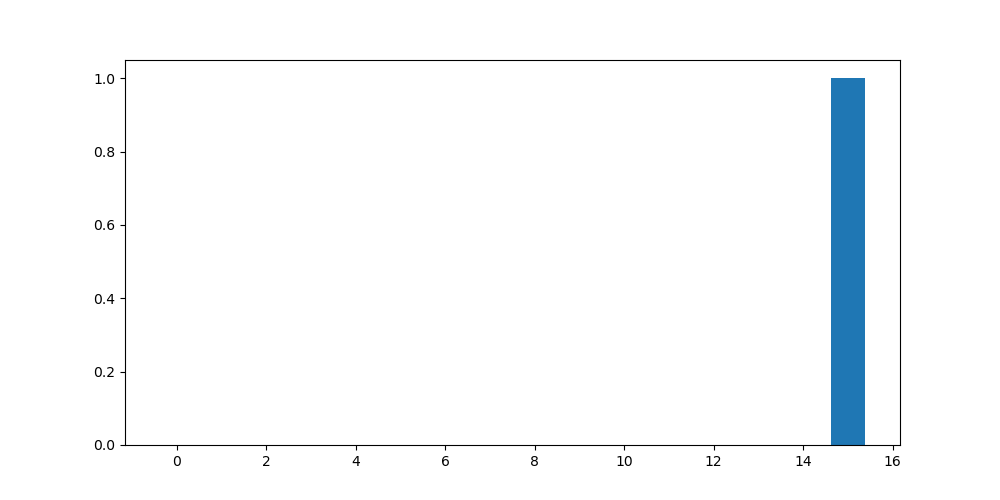}
\end{tabular}
\end{center}
\caption{Boltzmann-on-numbers distributions~\eqref{BoltzmannNumEqn} of
  the form $\boltzmannnum[16,7](i)$, with $K=7$ particles on the
  $N=16$ energy levels $\{0,1,\ldots,15\}$. The possible total
  energies $i\in\NNO$ must satisfy $0\leq i \leq 105$, where the
  highest value arises as $105 = (16-1)\cdot 7$.  The diagrams arise
  via steps of~$7$, with bar charts for $i=0,7,14,21$, on the top row,
  continuing downwards with $i=28,35,42,49$ etc. Notice that for
  energy total $i=0$ on the top-left the only possible multiset of
  particles is $7\ket{0}$, with all particles at energy level $0$,
  giving a singleton Boltzmann distribution $1\ket{0}$.  Similarly,
  for $i = 105$ at the lower-right, the only configuration is
  $7\ket{15}$, giving a singleton Boltzmann distribution
  $1\ket{15}$. }
\label{BoltzmannNumIterFig}
\end{figure}

As further illustration, several bar chart plots of the
Boltzmann-on-numbers distribution are given in
Figure~\ref{BoltzmannNumIterFig}.

%
%

We continue with basic properties of the Boltzmann-on-numbers
distribution. The first point describes the symmetry that is clearly
visible in Figure~\ref{BoltzmannNumIterFig}.

\begin{lemma}
\label{BoltzmannNumLem}
Fix numbers $N,K\geq 1$ and $0\leq i\leq (N\minnetje 1)\cdot K$.
\begin{enumerate}
\item \label{BoltzmannNumLemRevNum} The Boltzmann-on-numbers 
distribution is stable under reversal.
\[ \xymatrix@R-0.8pc@C+1pc{
\{0,1,\ldots,(N\minnetje 1)\cdot K\}\ar[d]_-{\revert_{(N-1)\cdot K}}^-{\cong}
   \ar[rr]^-{\boltzmannnum[N, K]} & &
   \Dst\big(\finset{N}\big)
   \ar[d]^-{\Dst(\revert_{N-1})}_-{\cong}
\\
\{0,1,\ldots,(N\minnetje 1)\cdot K\}\ar[rr]^-{\boltzmannnum[N, K]} & &
   \Dst\big(\finset{N}\big)
} \]

\item \label{BoltzmannNumLemMean} The mean of the Boltzmann-on-numbers 
distribution is $\frac{i}{K}$, independently of $N$:
\[ \begin{array}{rcl}
\mean\Big(\boltzmannnum[N, K](i)\Big)
& = &
\displaystyle\frac{i}{K}.
\end{array} \]
\end{enumerate}
\end{lemma}

\begin{myproof}
\begin{enumerate}
\item An abstract way to prove the statement, in the style of category
  theory, is to show first that the following `naturality' diagram
  commutes, and then use
  Lemma~\ref{BoltzmannMltLem}~\eqref{BoltzmannMltLemRev}.
\[ \vcenter{\xymatrix@R-0.8pc{
\Dst\Big(\natMlt[K](\finset{N})\Big)
   \ar[rr]^-{\flrn_*}\ar[d]_-{\Dst\natMlt(\revert_{N-1})}^-{\cong} 
   & &
   \Dst\big(\finset{N}\big)\ar[d]^-{\Dst(\revert_{N-1})}_-{\cong}
\\
\Dst\Big(\natMlt[K](\finset{N})\Big)\ar[rr]^-{\flrn_*} & &
   \Dst\big(\finset{N}\big)
}} \]

\noindent Alternatively, one can reason more concretely, like in the
proof of Lemma~\ref{BoltzmannMltLem}~\eqref{BoltzmannMltLemRev}:
\[ \begin{array}{rcl}
\lefteqn{\Dst(\revert_{N-1})\Big(\boltzmannnum[N, K](i)\Big)}
\\[+0.2em]
& \smash{\stackrel{\eqref{BoltzmannNumEqn}}{=}} &
\displaystyle\sum_{0\leq j < N} \, 
   \sum_{\varphi\in\natMlt[K](\finset{N}), \, \som(\varphi) = i} \,
   \frac{\coefm{\varphi}\cdot \varphi(j)}{C_{N}(K,i) \cdot K}
   \,\bigket{\revert_{N-1}(j)}
\\[+1em]
& = &
\displaystyle\sum_{0\leq j < N} \, 
   \sum_{\varphi\in\natMlt[K](\finset{N}), \, 
   \som(\natMlt(\revert_{N-1}(\varphi))) = i} \,
\\[+1.4em]
& & \hspace*{10em}\displaystyle
   \frac{\coefm{\natMlt(\revert_{N-1}(\varphi)}\cdot 
   \natMlt(\revert_{N-1}(\varphi))(j)}{C_{N}(K,i) \cdot K}
   \,\bigket{\revert_{N-1}(j)}
\\[+1em]
& = &
\displaystyle\sum_{0\leq j < N} \, 
   \sum_{\varphi\in\natMlt[K](\finset{N}), \, 
   \som(\varphi) = (N-1)\cdot K - i} \,
   \frac{\coefm{\varphi}\cdot 
   \varphi(\revert_{N-1}(j))}{K \cdot C_{N}(K,(N\minnetje 1)\cdot K - i)}
   \,\bigket{\revert_{N-1}(j)}
\\[+1em]
& = &
\displaystyle\sum_{0\leq j < N} \, 
   \sum_{\varphi\in\natMlt[K](\finset{N}), \, 
   \som(\varphi) = (N-1)\cdot K - i} \, \frac{\coefm{\varphi}\cdot 
   \varphi(j)}{K \cdot C_{N}(K,(N\minnetje 1)\cdot K - i)}
   \,\bigket{j}
\\[+1.2em]
& = &
\boltzmannnum[N, K]\Big((N\minnetje 1)\cdot K - i\Big).
\end{array} \]

\item We recall that the mean of a distribution
  $\omega\in\Dst(\finset{N})$ is defined as: $\mean(\omega) \coloneqq
  \sum_{j} \omega(j)\cdot j$.  This is applied to the
  Boltzmann-on-numbers distribution~\eqref{BoltzmannNumEqn}:
\[ \begin{array}[b]{rcl}
\mean\Big(\boltzmannnum[N, K](i)\Big)
& = &
\displaystyle\sum_{0\leq j < N} \, 
   \sum_{\varphi\in\natMlt[K,i](\finset{N})} \,
   \frac{\coefm{\varphi}\cdot \varphi(j)}{C_{N}(K,i) \cdot K}\cdot j
\\[+1.2em]
& = &
\displaystyle \sum_{\varphi\in\natMlt[K,i](\finset{N})} \,
   \frac{\coefm{\varphi}\cdot \sum_{0\leq j< N}\varphi(j) \cdot j}
   {C_{N}(K,i) \cdot K}
\\[+1.0em]
& = &
\displaystyle \sum_{\varphi\in\natMlt[K,i](\finset{N})} \,
   \frac{\coefm{\varphi}\cdot \som(\varphi)}{C_{N}(K,i) \cdot K}
\\[+1.0em]
& = &
\displaystyle \sum_{\varphi\in\natMlt[K,i](\finset{N})} \,
   \frac{\coefm{\varphi}\cdot i}{C_{N}(K,i) \cdot K}
\\[+1.4em]
& \smash{\stackrel{\eqref{NnomialMltEqn}}{=}} &
\displaystyle \frac{C_{N}(K,i) \cdot i}{C_{N}(K,i) \cdot K}
\hspace*{\arraycolsep}=\hspace*{\arraycolsep}
\frac{i}{K}.
\end{array} \eqno{\square} \]
\end{enumerate}
\end{myproof}

We conclude this section with two alternative formulations of the
Boltzmann-on-numbers distributions. The second point below show that
this distribution can be obtained directly from the uniform
distribution on sequences (microstates), without considering any
multisets (macrostates).

\begin{proposition}
\label{BoltzmannNumAltProp}
let $N\geq 1$ and $K\geq 0$ with $0\leq i \leq (N\minnetje 1)\cdot K$.
\begin{enumerate}
\item \label{BoltzmannNumAltPropNom} The Boltzmann-on-numbers
  distribution can be expressed purely in terms of $N$-nomials:
\begin{equation}
\label{BoltzmannNumAltPropNomEqn}
\begin{array}{rcl}
\boltzmannnum[N,K](i)
& = &
\displaystyle\sum_{0\leq j < \min(N,i+1)} \, \frac{C_{N}(K-1, i-j)}{C_{N}(K,i)}
   \,\bigket{j}.
\end{array}
\end{equation}

\item \label{BoltzmannNumAltPropProj} Recall the set $S\subseteq
  \finset{N}^{K}$ with the uniform distribution distribution
  $\unif\in\Dst(S)$ from
  Lemma~\ref{BoltzmannMltLem}~\eqref{BoltzmannMltLemIm}. For any $0
  \leq n < N$ let $\pi_{n} \colon \finset{N}^{K} \rightarrow
  \finset{N}$ be the $n$-th projection from a sequence: $\pi_{n}(v) =
  v_{n}$. Then, for each $n$,
\[ \begin{array}{rcl}
\Dst\big(\pi_{n}\big)(\unif)
& = &
\boltzmannnum[N,K](i).
\end{array} \]


\ignore{

N = 4
Nsp = range_sp(N)
i = 6
for K in range(2, 10):
    bom = boltzmann_on_multisets(N,K,frac=True)(i)
    print("--> ", K)
    print( bom )
    avg = empty_multiset(Nsp, frac=True)
    for m in bom.support_list():
        avg += bom(m) * DState(m.matrix, Nsp, frac=True)
    print( Frac(1,K) * avg )
    print( boltzmann_on_numbers(N,K,frac=True)(i) )

}
\end{enumerate}
\end{proposition}

\begin{myproof}
\begin{enumerate}
\item The result follows from the equation:
\[ \begin{array}[b]{rcl}
\displaystyle\sum_{\varphi\in\natMlt[K,i](\finset{N})} \,
   \coefm{\varphi}\cdot \varphi(j)
& \smash{\stackrel{\eqref{MltCoefEqn}}{=}} &
\displaystyle \sum_{\varphi\in\natMlt[K,i](\finset{N}), \,
   \varphi(j) \neq 0} \,
   \frac{K!}{\varphi(j)!\cdot \prod_{m\neq j} \varphi(m)!} \cdot \varphi(j)
\\[+1em]
& = &
\displaystyle K \cdot 
   \sum_{\varphi\in\natMlt[K,i](\finset{N}), \, \varphi(j) \neq 0} \,
   \frac{(K-1)!}{(\varphi(j)-1)!\cdot \prod_{m\neq j} \varphi(m)!}
\\[+1.2em]
& = &
\displaystyle K \cdot 
   \sum_{\varphi\in\natMlt[K](\finset{N}), \, \som(\varphi) = i, \,
   \varphi(j) \neq 0} \, \coefm{\varphi - 1\ket{j}}
\\[+1.2em]
& = &
\displaystyle K \cdot 
   \sum_{\psi\in\natMlt[K-1](\finset{N}), \, \som(\psi) = i-j} \,
   \coefm{\psi}
\\[+1.4em]
& \smash{\stackrel{\eqref{NnomialMltEqn}}{=}} &
K \cdot C_{N}\big(K-1, i-j\big).
\end{array} \]

\item We use the previous point:
\[ \begin{array}[b]{rcl}
\Dst\big(\pi_{n}\big)(\unif)
& = &
\displaystyle\sum_{v\in S} \, \frac{1}{C_{N}(K,i)} \, \bigket{\pi_{n}(v)}
\\
& = &
\displaystyle \sum_{0\leq j < N} \,
   \sum_{v\in S, \, \pi_{n}(v) = j} \, \frac{1}{C_{N}(K,i)} \, \bigket{j}
\\
& = &
\displaystyle \sum_{0 \leq j < \min(N, i+1)} \,
   \sum_{w\in \finset{N}^{K-1}, \, \som(w) = i-j} \, \frac{1}{C_{N}(K,i)} \, \bigket{j}
\\[+1.4em]
& \smash{\stackrel{\eqref{NnomialSeqEqn}}{=}} &
\displaystyle\sum_{0\leq j < \min(N,i+1)} \, \frac{C_{N}(K-1, i-j)}{C_{N}(K,i)}
   \,\bigket{j}.
\\
& \smash{\stackrel{\eqref{BoltzmannNumAltPropNomEqn}}{=}} &
\boltzmannnum[N,K](i).
\end{array} \eqno{\square} \]

\end{enumerate}
\end{myproof}

\ignore{

E = 15
N = 11
K = 4
S = [ (n1, n2, n3, n4) for n1 in range(N) for n2 in range(N)
      for n3 in range(N) for n4 in range(N) if n1 + n2 + n3 + n4 == E ]

SP = Space(S)

u = uniform_distribution(SP, frac=True)

print("")
b = 2
print( u >= DPred(lambda ns: 1 if ns[1] == b else 0, SP, frac=True) )
print( Frac(nomial(N, K-1, E-b), nomial(N, K, E)) )

print("\nAs distribution")
print( DState(lambda i: u >= DPred(lambda ns: 1 if ns[1] == i else 0, SP, frac=True), range_sp(N), frac=True) )

print("")

print( DState(lambda i: Frac(nomial(N, K-1, E-i), nomial(N, K, E)), range_sp(N), frac=True) )

print("") 

print( boltzmann_on_numbers_rec(N, K, frac=True)(E) )

}

\section{The Boltzmann distribution on energy levels}\label{BoltzmannEneSec}

As we briefly discussed above, physicists see the Boltzmann-on-numbers
distribution $\boltzmannnum[N, K](i)$ as a distribution on energy
levels $\{0,\ldots,N\minnetje 1\} = \finset{N}$, arising from $K$
particles with total energy $i$. The pictures that they look at
typically have shapes as in the top row of
Figure~\ref{BoltzmannNumIterFig}. They see the number $N$ of energy
levels as a number that can go to infinity, in order to get an
approximation of a continuous distribution, see
Subsection~\ref{BoltzmannApproxSubsec} for an impression. Here the
focus is not such approximation, but on the precise form of the
relevant discrete probability distributions.  So far we used a
constraint on the number $i$ of the form $0 \leq i \leq (N\minnetje
1)\cdot K$. It turns out that the constructions that we used so far
can be simplified when we assume $i < N$. This is a crucial
observation, elaborated in the theorem below.  The assumption $i < N$
makes sense, when $N$ is seen as a number that can become arbitrarily
big.

\begin{theorem}
\label{BelowThm}
Consider numbers $N,K \geq 1$.
\begin{enumerate}
\item \label{BelowThmSum} For $0 \leq n \leq N$,
\[ \begin{array}{rcl}
\displaystyle \sum_{0\leq i < n} \, C_{N}(K,i)
& = &
\displaystyle\frac{n}{K}\cdot\bibinom{K}{n}.
\end{array} \]

\noindent The right-hand-side does not depend on $N$, intuitively because
if $\som(\vec{m})$ equals $i$, for $\vec{m} \in \finset{N}^{K}$ where
$i \leq n < N$, then the number $N$ is big enough to so that the
requirement $m_{j} \in \finset{N}$ does not impose any restriction.

\item \label{BelowThmNomial} When $i < N$, there is a simple formula
for the $N$-nomial, namely:
\[ \begin{array}{rcl}
C_{N}(K,i)
& = &
\displaystyle\bibinom{K}{i}.
\end{array} \]
\end{enumerate}
\end{theorem}

\begin{myproof}
\begin{enumerate}
\item By induction on $K\geq 1$. First,
\[ \begin{array}{rcl}
\displaystyle \sum_{0\leq i\leq n} \, C_{N}(1,i)
\hspace*{\arraycolsep}\smash{\stackrel{\eqref{NnomialSeqEqn}}{=}}\hspace*{\arraycolsep}
\sum_{0\leq i\leq n} \, 
   \Bigsetsize{\bigsetin{\vec{n}}{\finset{N}^{1}}{\som(\vec{n}) = i}}
& = &
\displaystyle \sum_{0\leq i\leq n} \, 1
\\
& = &
n
\hspace*{\arraycolsep}=\hspace*{\arraycolsep}
\displaystyle\frac{n}{1}\cdot\bibinom{1}{n}.
\end{array} \]

\noindent Next we use that we sum over $i < n \leq N$ in the first
line below. When we count sequences of length $K\plusje 1$ with sum $i
< n$, then all elements of the sequence fit in the set $\finset{N} =
\{0,\ldots,N\minnetje 1\}$; hence we can express the number
$C_{N}(K\plusje 1, i)$ of sequences of length $K\plusje 1$ as sum of
numbers of sequences of length $K$.
\[ \begin{array}{rcl}
\displaystyle \sum_{0\leq i < n} \, C_{N}\big(K+1,i\big)
& = &
\displaystyle \sum_{0\leq i < n} \, \sum_{0\leq j \leq i} C_{N}\big(K,i-j\big)
\\[+1em]
& = &
\displaystyle \sum_{0\leq i < n} \, \sum_{0\leq j \leq i} C_{N}\big(K,j\big)
\\[+1.2em]
& \smash{\stackrel{\textrm{(IH)}}{=}} &
\displaystyle \displaystyle \sum_{0\leq i < n} \, 
   \frac{i+1}{K}\cdot\bibinom{K}{i+1}
\\[+1em]
& = &
\displaystyle \displaystyle \sum_{0\leq i < n} \, 
   \frac{i+1}{K}\cdot\frac{(K+i)!}{(K-1)!\cdot (i+1)!}
\\[+1em]
& = &
\displaystyle \displaystyle \sum_{1\leq i < n} \, \bibinom{K+1}{i}
\\[+1.2em]
& = &
\displaystyle\bibinom{n}{K+1} \qquad 
   \mbox{by Lemma~\ref{ChooseSumLem}~\eqref{ChooseSumLemZero}}
\\[+0.8em]
& = &
\displaystyle\frac{(K+n)!}{(n-1)!\cdot (K+1)!}
\hspace*{\arraycolsep}=\hspace*{\arraycolsep}
\frac{n}{K+1}\cdot\bibinom{K+1}{n}.
\end{array} \]

\item By the previous item, for $i < N$, so that $i\plusje 1 \leq N$,
\[ \hspace*{-0.6em}\begin{array}[b]{rcl}
C_{N}(K,i)
& = &
\displaystyle \left(\sum_{0\leq j < i+1} \, C_{N}\big(K,j\big)\right)
  - \left(\sum_{0\leq j < i} \, C_{N}\big(K,j\big)\right)
\\[+1.4em]
& = &
\displaystyle\frac{i+1}{K}\cdot\bibinom{K}{i+1} -
   \frac{i}{K}\cdot\bibinom{K}{i}
\\[+1em]
& = &
\displaystyle\frac{(K+i)!}{K!\cdot i!} - \frac{(K+i-1)!}{K!\cdot (i-1)!}
\\[+1em]
& = &
\displaystyle\frac{(K+i-1)!}{K!\cdot (i-1)!} \cdot 
   \left(\frac{K+i}{i} - 1\right)
\hspace*{\arraycolsep}=\hspace*{\arraycolsep}
\displaystyle\frac{(K+i-1)!}{(K-1)!\cdot i!}
\hspace*{\arraycolsep}=\hspace*{\arraycolsep}
\displaystyle\bibinom{K}{i}\!.
\end{array} \eqno{\square} \]
\end{enumerate}
\end{myproof}

We finally define the Boltzmann-on-energy distribution, as a a
distribution on energy levels, via the simplifications from the
previous theorem. To emphasise that we consider the number $i$ as the
total energy we shall now rename it to $e$.

\begin{definition}
\label{BoltzmannEneDef}
Let natural numbers $E \geq 1$ and $K \geq 2$ be given, where $K$
stands for the number of particles and $E$ for the maximal / total
energy. We define the Boltzmann-on-energy function $\boltzmannene[E]
\colon \setin{K}{\NNO}{K\geq 2} \rightarrow \Dst\big(\{0,\ldots,
E\}\big)$ as:
\begin{equation}
\label{BoltzmannEneEqn}
\begin{array}{rcl}
\boltzmannene[E](K)
\hspace*{\arraycolsep}\coloneqq\hspace*{\arraycolsep}
\boltzmannnum\big[E\plusje 1, K\big](E)
& = &
\displaystyle\sum_{0\leq j \leq E} \, 
      \frac{\big(\binom{K-1}{E-j}\big)}{\big(\binom{K}{E}\big)}\,\bigket{j}.
\end{array}
\end{equation}
\end{definition}

Notice that in the Boltzmann-on-numbers distribution in the middle
of~\eqref{BoltzmannEneEqn} we use $N = E+1$ and $i = E$, so that
condition $i < N$ from Theorem~\ref{BelowThm} holds by construction.
This justifies the simple form at the right-hand-side
of~\eqref{BoltzmannEneEqn}, via
Proposition~\ref{BoltzmannNumAltProp}~\eqref{BoltzmannNumAltPropNom}
and Theorem~\ref{BelowThm}~\eqref{BelowThmNomial}:
\[ \begin{array}[b]{rcl}
\boltzmannene[E](K)
\hspace*{\arraycolsep}=\hspace*{\arraycolsep}
\boltzmannnum\big[E\plusje 1, K\big](E)
& \smash{\stackrel{\eqref{BoltzmannNumAltPropNomEqn}}{=}} &
\displaystyle\sum_{0\leq j < \min(E+1,E+1)} \, \frac{C_{N}(K-1, E-j)}{C_{N}(K,E)}
   \,\bigket{j}.
\\[+1.2em]
& = &
\displaystyle\sum_{0\leq j \leq E} \, 
      \frac{\big(\binom{K-1}{E-j}\big)}{\big(\binom{K}{E}\big)}\,\bigket{j}.
\end{array} \]

\begin{figure}
\begin{center}
\hspace*{-1em}\begin{tabular}{ccc}
\includegraphics[scale=0.18]{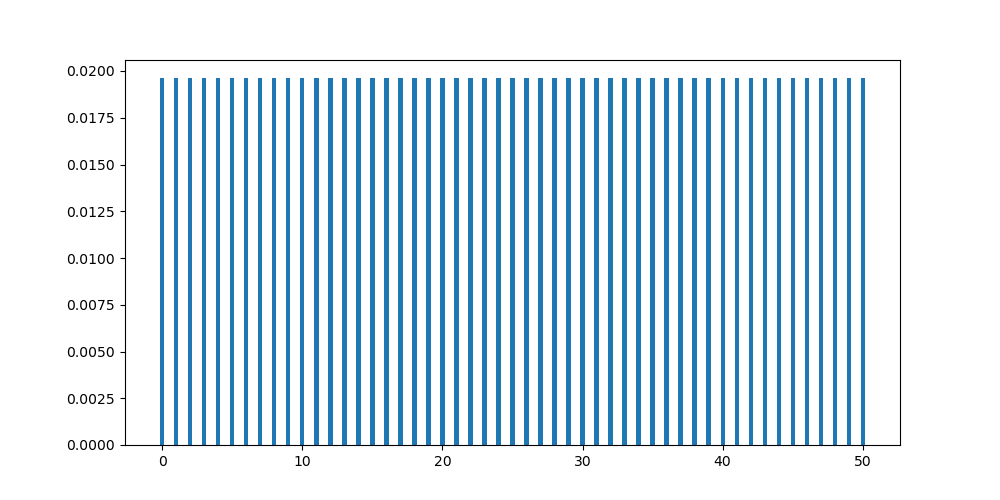}
   \hspace*{-1em} & \hspace*{-1em}
\includegraphics[scale=0.18]{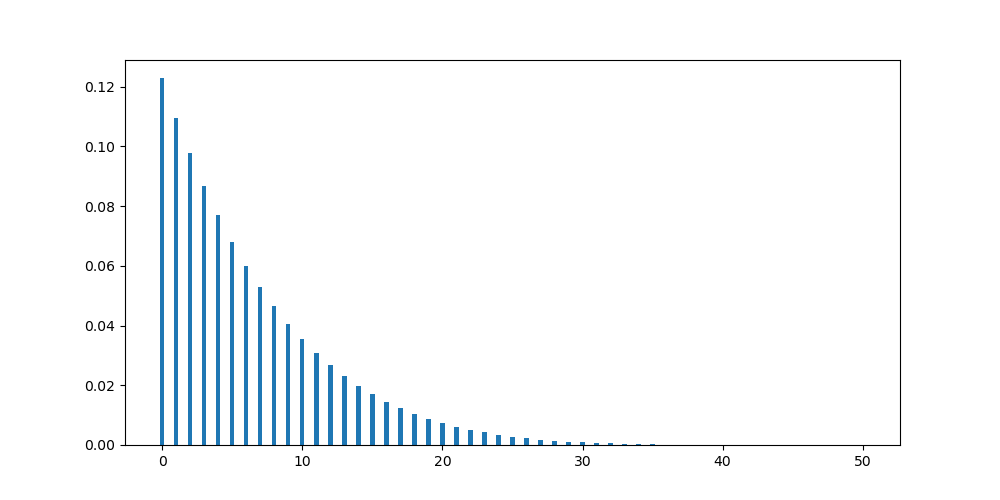}
   \hspace*{-1em}  & \hspace*{-1em}
\includegraphics[scale=0.18]{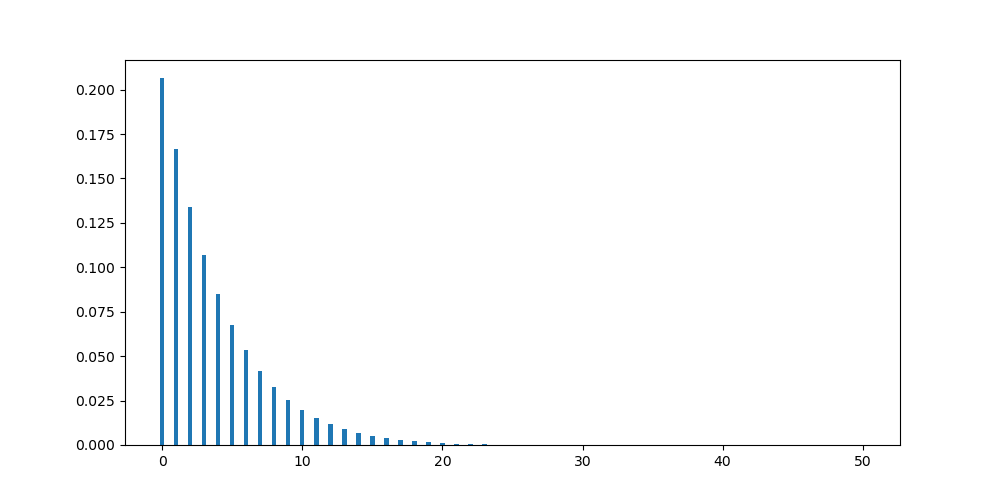}
\\
\includegraphics[scale=0.18]{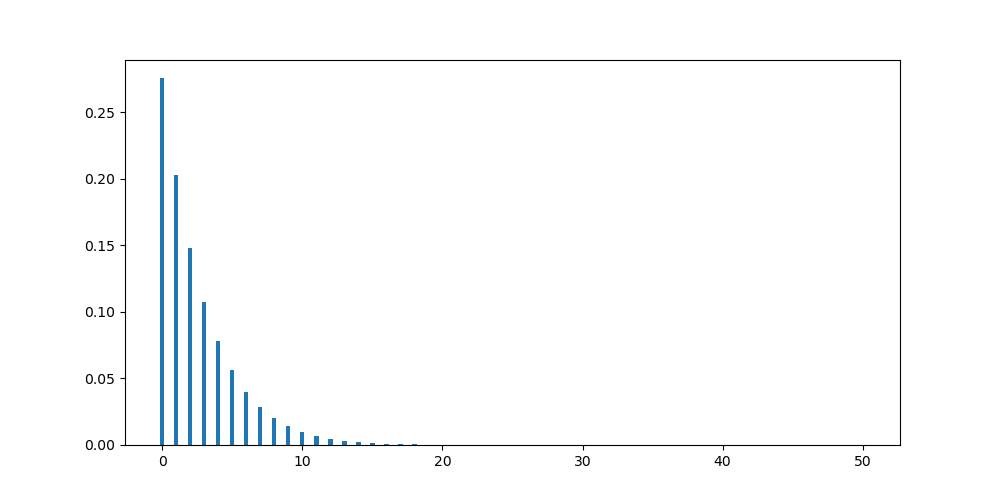}
   \hspace*{-1em} & \hspace*{-1em}
\includegraphics[scale=0.18]{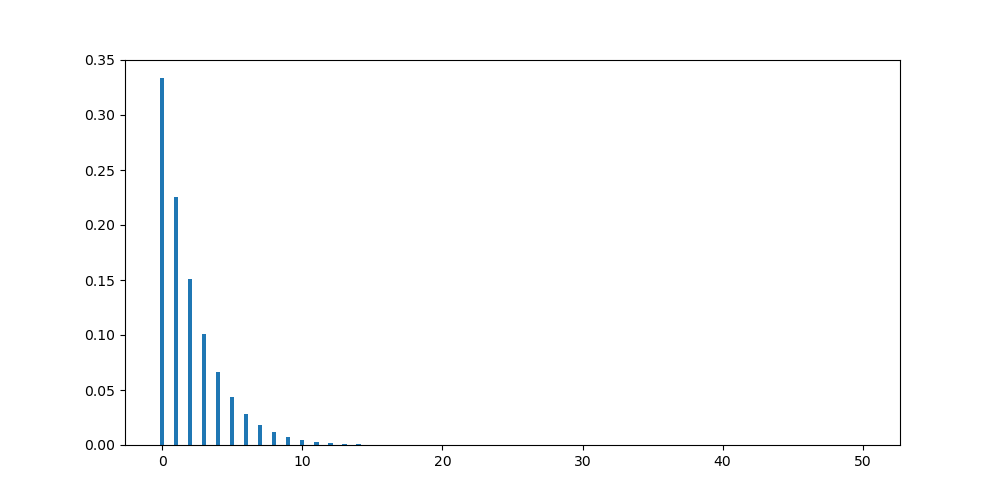}
   \hspace*{-1em}  & \hspace*{-1em}
\includegraphics[scale=0.18]{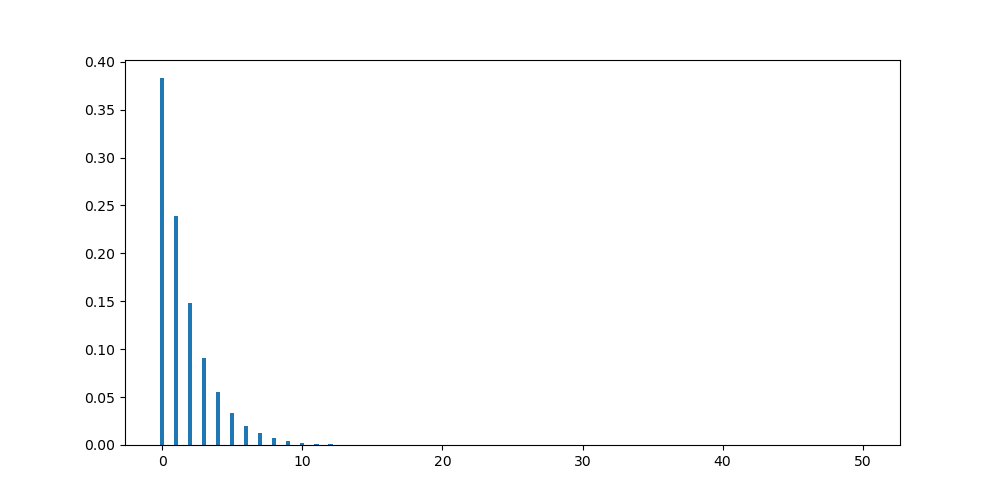}
\\
\includegraphics[scale=0.18]{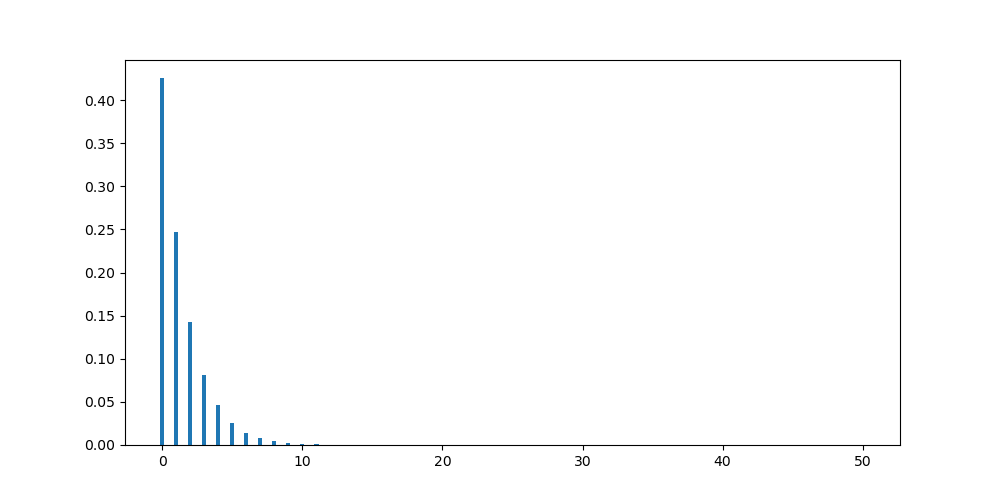}
   \hspace*{-1em} & \hspace*{-1em}
\includegraphics[scale=0.18]{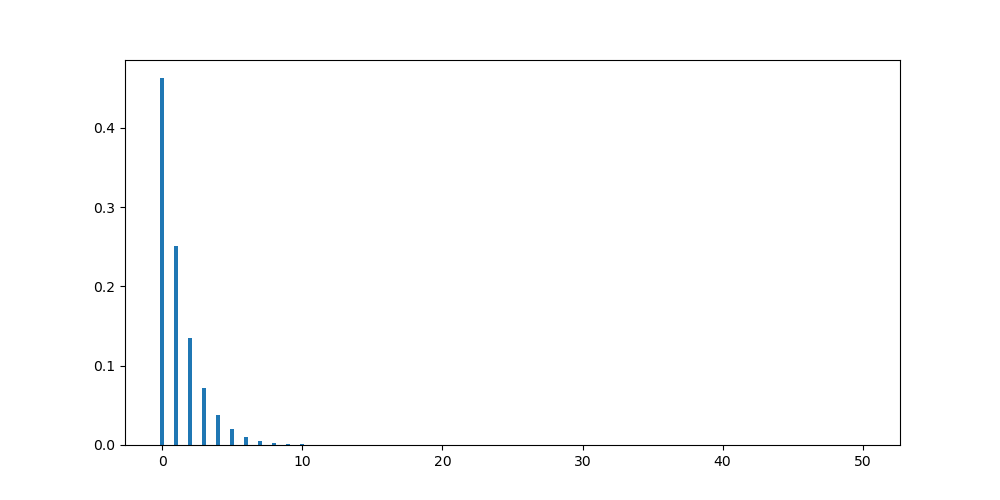}
   \hspace*{-1em}  & \hspace*{-1em}
\includegraphics[scale=0.18]{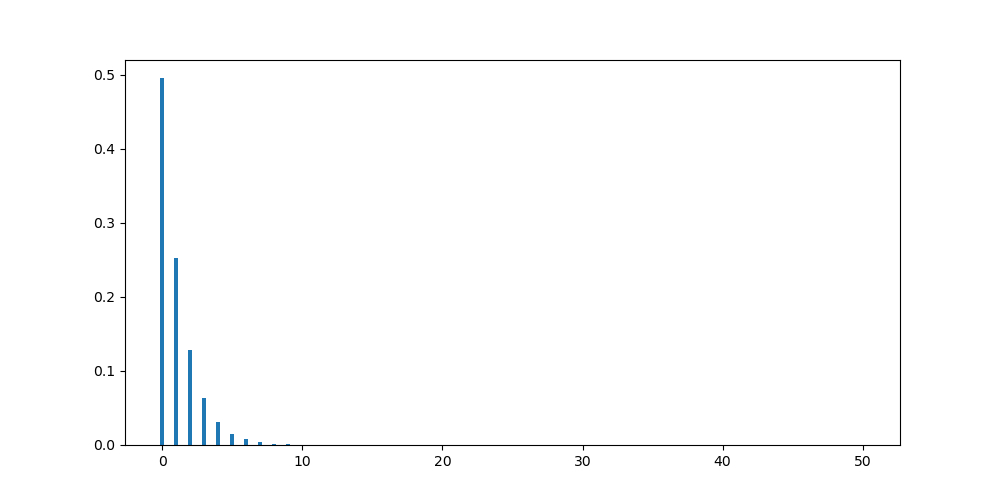}
\end{tabular}
\end{center}
\caption{Boltzmann-on-energy distributions~\eqref{BoltzmannEneEqn} of
  the form $\boltzmannene[50](K)$ on energy levels
  $\{0,1,\ldots,50\}$, with maximal/total energy $E=50$ for $K$
  particles. The top row shows the bar charts for $K=2,8,14$, the
  middle row for $K=20,26,32$ and the bottom row for $K=38,44,50$. We
  see that the higher the number of particles, the more they
  concentrate on the lower energy levels --- since the sum of all all
  their energy levels must remain equal to~$E=50$. For $K=2$ the
  distribution is uniform, see
  Proposition~\ref{BoltzmannEneProp}~\eqref{BoltzmannEnePropUnif} for
  the general case.}
\label{BoltzmannEneIterFig}
\end{figure}

\noindent Figure~\ref{BoltzmannEneIterFig} gives a systematic
description of these Boltzmann-on-energy distributions. But we also
(re-)elaborate our running example.

\begin{example}
\label{PhysicsEneEx}
The Boltzmann-on-integers distribution $\boltzmannnum[4, 4](3)$,
computed (twice) in Example~\ref{PhysicsNumEx}, can be computed a
third time, but now as Boltzmann-on-energy:
\[ \begin{array}{rcl}
\boltzmannene[3](4)
\hspace*{\arraycolsep}=\hspace*{\arraycolsep}
\boltzmannnum[4, 4](3)
& \smash{\stackrel{\eqref{BoltzmannEneEqn}}{=}} &
\displaystyle\sum_{0\leq j \leq 3} \, 
      \frac{\big(\binom{3}{3-j}\big)}{\big(\binom{4}{3}\big)}\,\bigket{j}
\\[+1.2em]
& = &
\displaystyle\frac{\big(\binom{3}{3}\big)}{20}\,\bigket{0} + 
   \frac{\big(\binom{3}{2}\big)}{20}\,\bigket{1} + 
   \frac{\big(\binom{3}{1}\big)}{20}\,\bigket{2} + 
   \frac{\big(\binom{3}{0}\big)}{20}\,\bigket{3} 
\\[+0.4em]
& = &
\frac{10}{20}\bigket{0} + \frac{6}{20}\bigket{1} + 
   \frac{3}{20}\bigket{2} + \frac{1}{20}\bigket{3}
\\[+0.2em]
& = &
\frac{1}{2}\bigket{0} + \frac{3}{10}\bigket{1} + 
   \frac{3}{20}\bigket{2} + \frac{1}{20}\bigket{3}.
\end{array} \]
\end{example}


\begin{proposition}
\label{BoltzmannEneProp}
Let $E\geq 1$ and $K \geq 2$.
\begin{enumerate}
\item \label{BoltzmannEnePropUnif} For $K=2$ particles the
  Boltzmann-on-energy distribution is uniform, of the form:
\[ \begin{array}{rcl}
\boltzmannene[E](2)
& = &
\displaystyle\sum_{0\leq j\leq E} \, \frac{1}{E+1} \, \bigket{j}.
\end{array} \]

\item \label{BoltzmannEnePropMean} The mean of the Boltzmann-on-energy
  distribution is given by:
\[ \begin{array}{rcl}
\mean\Big(\boltzmannene[E](K)\Big)
& = &
\displaystyle\frac{E}{K}.
\end{array} \]

\item \label{BoltzmannEnePropVar} And its variance is:
\[ \begin{array}{rcl}
\Var\Big(\boltzmannene[E](K)\Big)
& = &
\displaystyle\frac{E\cdot (E+K)\cdot (K-1)}{K^{2}\cdot (K+1)}.
\end{array} \]
\end{enumerate}
\end{proposition}

\begin{myproof}
\begin{enumerate}
\item For each $0\leq j\leq e$ one has:
\[ \begin{array}{rcccccl}
\boltzmannene[E](2)(j)
& = &
\displaystyle\frac{\big(\binom{2-1}{E-j}\big)}{\big(\binom{2}{E}\big)}
& = &
\displaystyle\frac{\frac{(E-j)!}{0!\cdot (E-j)!}}
   {\frac{(e+1)!}{1!\cdot e!}}
& = &
\displaystyle\frac{1}{E+1}.
\end{array} \]

\item The equation about the mean is obtained from
  Lemma~\ref{BoltzmannNumLem}~\eqref{BoltzmannNumLemMean}, since:
\[ \begin{array}{rcccl}
\mean\Big(\boltzmannene[E](K)\Big)
& \smash{\stackrel{\eqref{BoltzmannEneEqn}}{=}} &
\mean\Big(\boltzmannnum\big[E\plusje 1, K\big](E)\Big)
& = &
\displaystyle\frac{E}{K}.
\end{array} \]

\item Proving the variance equation is more work. We first
show that:
\[ \begin{array}{rcl}
\displaystyle\sum_{0\leq j \leq E} \, 
   \frac{\big(\binom{K-1}{E-j}\big)}{\big(\binom{K}{E}\big)} \cdot j^{2}
& = &
\displaystyle\frac{E\cdot (K-1+2E)}{K\cdot (K+1)}.
\end{array} \eqno{(*)} \]

\noindent Indeed:
\[ \hspace*{-0.5em}\begin{array}{rcl}
\lefteqn{\displaystyle\sum_{0\leq j \leq E} \, 
   \frac{\big(\binom{K-1}{E-j}\big)}{\big(\binom{K}{E}\big)} \cdot j^{2}}
\\[+0.6em]
& = &
\displaystyle\frac{1}{\big(\binom{K}{E}\big)} \cdot
   \sum_{0\leq j \leq E} \, \bibinom{K-1}{j} \cdot (E-j)^{2}
\\[+1.2em]
& = &
\displaystyle\frac{1}{\big(\binom{K}{E}\big)} \cdot \left[\,
   E^{2} \cdot \sum_{0\leq j \leq E} \, \bibinom{K-1}{j} -
   2E \cdot \sum_{0\leq j \leq E} \, \bibinom{K-1}{j} \cdot j \right.
\\
& & \hspace*{8em} + \; \displaystyle \left.
   \sum_{0\leq j \leq E} \, \bibinom{K-1}{j} \cdot j^{2} \, \right]
\\[+1.2em]
& = &
\displaystyle\frac{1}{\big(\binom{K}{E}\big)} \cdot \left[\,
   E^{2} \cdot \bibinom{E+1}{K-1} 
   - 2E\cdot (K\minnetje 1) \cdot \bibinom{E}{K} \right.
\\
& & \hspace*{8em} + \; \displaystyle \left.
   K\cdot (K\minnetje 1) \cdot \bibinom{E-1}{K+1} +
   (K\minnetje 1)\cdot \bibinom{E}{K} \, \right]
\\[+1.4em]
& = &
\displaystyle\frac{(K\minnetje 1)!\cdot E!}{(K\plusje E \minnetje 1)!} \cdot
   \left[\, E^{2}\cdot\frac{(K\plusje E \minnetje 1)!}
      {E!\cdot (K\minnetje 1)!}
   - (2E\minnetje 1)\cdot (K\minnetje 1) \cdot 
     \frac{(K\plusje E \minnetje 1)!}{(E\minnetje 1)!\cdot K!} \right.
\\[+0.8em]
& & \hspace*{8em} \displaystyle \left.
   +\; K\cdot (K\minnetje 1) \cdot 
     \frac{(K\plusje E \minnetje 1)!}{(E\minnetje 2)!\cdot (K\plusje 1)!}
   \,\right]
\\[+0.6em]
& = &
\displaystyle E^{2} 
   - \frac{E\cdot (2E \minnetje 1) \cdot (K \minnetje 1)}{K}
   + \frac{(K\minnetje 1)\cdot E \cdot (E \minnetje 1)}{K + 1}
\\[+0.8em]
& = &
\displaystyle\frac{E^{2}\cdot K^{2} \plusje E^{2}\cdot K \minnetje 2E\cdot K^{2}
   \plusje 2E^{2} \plusje E\cdot K^{2} \minnetje E \plusje E^{2} \cdot K^{2} 
   \minnetje E \cdot K^{2} \minnetje E^{2} \cdot K \plusje E \cdot K}
   {K\cdot (K+1)}
\\[+0.8em]
& = &
\displaystyle\frac{2E^{2} - E + E\cdot K}{K\cdot (K+1)}
\hspace*{\arraycolsep}=\hspace*{\arraycolsep}
\frac{E\cdot (K-1+2E)}{K\cdot (K+1)}.
\end{array} \]

\noindent Now we come to variance itself, using the description of the mean:
\[ \begin{array}[b]{rcl}
\Var\Big(\boltzmannene[E](K)\Big)
& = &
\displaystyle\sum_{0\leq j \leq E} \, 
   \frac{\big(\binom{K-1}{E-j}\big)}{\big(\binom{K}{E}\big)} \cdot j^{2} \;-\;
   \left(\sum_{0\leq j \leq E} \, 
   \frac{\big(\binom{K-1}{E-j}\big)}{\big(\binom{K}{E}\big)} \cdot j\right)^{2}
\\[+1.2em]
& \smash{\stackrel{(*)}{=}} &
\displaystyle\frac{E\cdot (K-1+2E)}{K\cdot (K+1)} - 
   \left(\frac{E}{K}\right)^{2}
\\[+0.8em]
& = &
\displaystyle\frac{e\cdot (K-1+2e) \cdot K - e^{2}\cdot (K+1)}{K^{2}\cdot (K+1)}
\\[+0.8em]
& = &
\displaystyle\frac{E\cdot (E+K)\cdot (K-1)}{K^{2}\cdot (K+1)}.
\end{array} \eqno{\square} \]
\end{enumerate}
\end{myproof}

\section{Final observations}\label{ObservationSec}

Having seen precise derivations of various discrete Boltzmann
distributions we lean back for some broader observations.  The style
of section is a bit less thorough and part of it may be elaborated
further in future publications.

\subsection{Approximations of the discrete Boltzmann 
   distributions}\label{BoltzmannApproxSubsec}

What is called the discrete or continuous Boltzmann distribution in
statistical mechanics is an approximation of what we have described
as the Boltzmann-on-energy distribution. We briefly sketch how this
can work, in several ways.

For total energy $E$ with $K$ particles the Boltzmann-on-energy
distribution $\boltzmannene[E](K)$ has mean $\mu \coloneqq \frac{E}{K}$,
see Proposition~\ref{BoltzmannEneProp}~\eqref{BoltzmannEnePropMean}.
We can look at the ratio between two successive probabilities, at
$j$ and at $j\plusje 1$. This gives, when we substitute $E = K\cdot \mu$,
\[ \begin{array}{rcl}
\displaystyle\frac{\boltzmannene[E](K)(j\plusje 1)}{\boltzmannene[E](K)(j)}
\hspace*{\arraycolsep}\smash{\stackrel{\eqref{BoltzmannEneEqn}}{=}}\hspace*{\arraycolsep}
\displaystyle\frac{\big(\binom{K-1}{E-(j+1)}\big)}{\big(\binom{K-1}{E-j}\big)}
& = &
\displaystyle\frac{(K+E-j-3)! \cdot (E-j)!}{(K+E-j-2)! \cdot (E-j-1)!}
\\[+1em]
& = &
\displaystyle\frac{E-j}{K+E-j-2}
\\[+0.8em]
& = &
\displaystyle\frac{K\cdot \mu-j}{K+K\cdot \mu-j-2}
  \;\xrightarrow{\quad K\rightarrow\infty\quad} \frac{\mu}{\mu+1}.
\end{array} \]

\noindent This suggests that we can use as approximation of the
Boltzmann-on-energy distribution, in presence of many particles, the
normalisation of:
\begin{equation}
\label{DiscIterApprox}
\sum_{0\leq j \leq E} \, \left(\frac{\mu}{\mu+1}\right)^{j}\,\bigket{j}
\end{equation}

\noindent Distributions with a constant ratio between successive
probabilities typically occur for (normalisations of) distributions
of the form $\sum_{j} e^{-b\cdot j} \ket{j}$, for a suitable constant
$b$. Indeed, the associated ratio is:
\[ \begin{array}{rcccl}
\displaystyle\frac{e^{-b\cdot (j+1)}}{e^{-b\cdot j}}
& = &
e^{-b\cdot (j+1) + b\cdot j}
& = &
e^{-b}.
\end{array} \]

\noindent In our Boltzmann-on-energy situation we can solve $e^{-b} =
\frac{\mu}{\mu + 1}$ and get:
\[ \begin{array}{rcccccccl}
b
& = &
\displaystyle -\ln\left(\frac{\mu}{\mu+1}\right)
& = &
\displaystyle \ln\left(\frac{\mu+1}{\mu}\right)
& = &
\displaystyle \ln\left(1 + \frac{1}{\mu}\right)
& \approx &
\displaystyle\frac{1}{\mu}.
\end{array} \]

\noindent The latter approximation is common for small
$\frac{1}{\mu}$, that is, for large $\mu$. We thus get as alternative
discrete approximation of the Boltzmann-on-energy distribution, the
normalisation of:
\begin{equation}
\label{DiscExpApprox}
\sum_{0\leq j \leq E} \, e^{-\frac{j}{\mu}}\,\bigket{j}
\end{equation}

The principle of maximum entropy introduced in~\cite{Jaynes57} is used
in statistical mechanics to derive the equilibrium distributions. We
sketch how this works in the current situation, with $N$ energy levels
and a given mean $\mu$. The aim is to find the distribution $\omega$
on $\{0,\ldots,E\}$ with $\mean(\omega) = \mu$ and maximum entropty
$-\sum_{i} \omega(i) \cdot \ln\big(\omega(i)\big)$. This can be done
via Lagrange's multiplier method, as described for a dice
in~\cite[Example~5.3]{DillB10}. We give an impression of how this
works.

Take as function:
\[ \begin{array}{rcl}
H\big(\vec{r}, \alpha, \beta\big)
& \coloneqq &
\displaystyle-\sum_{0\leq j \leq E} r_{j}\cdot\ln(r_{j}) - 
   \alpha\cdot\left(1 - \sum_{0\leq j \leq E} r_{j}\right) -
   \beta\cdot\left(\mu - \sum_{0\leq j \leq E} r_{j}\cdot j\right).
\end{array} \]

\auxproof{
The partial derivatives are:
\[ \begin{array}{rcl}
\displaystyle\frac{\partial H}{\partial r_{j}}\big(\vec{r}, \alpha, \beta\big)
& = &
-\ln(r_{j}) - 1 + \alpha + \beta\cdot j
\\[+1em]
\displaystyle\frac{\partial H}{\partial \alpha}\big(\vec{r}, \alpha, \beta\big)
& = &
\displaystyle 1 - \sum_{0\leq j \leq E} r_{j}
\\[+1em]
\displaystyle\frac{\partial H}{\partial \beta}\big(\vec{r}, \alpha, \beta\big)
& = &
\displaystyle 1 - \sum_{0\leq j \leq E} r_{j}\cdot j.
\end{array} \]
}

\noindent Taking partial derivatives and setting them to zero yields,
when we abbreviate $x = e^{\beta}$.
\[ \begin{array}{rcl}
1
& = &
\displaystyle \sum_{0\leq j \leq E}\, r_{j}
\hspace*{\arraycolsep}=\hspace*{\arraycolsep}
\displaystyle \sum_{0\leq j \leq E}\, e^{-1 + \alpha + \beta\cdot j}
\hspace*{\arraycolsep}=\hspace*{\arraycolsep}
\displaystyle \sum_{0\leq j \leq E}\, \frac{(e^{\beta})^{i}}{e^{1-\alpha}}
\hspace*{\arraycolsep}=\hspace*{\arraycolsep}
\displaystyle \sum_{0\leq j \leq E}\, \frac{x^{i}}{e^{1-\alpha}}
\\[+1em]
\mu
& = &
\displaystyle \sum_{0\leq j \leq E}\, r_{j}\cdot j
\hspace*{\arraycolsep}=\hspace*{\arraycolsep}
\displaystyle \sum_{0\leq j \leq E}\, e^{-1 + \alpha + \beta\cdot j}\cdot j
\hspace*{\arraycolsep}=\hspace*{\arraycolsep}
\displaystyle \sum_{0\leq j \leq E}\, \frac{(e^{\beta})^{j}\cdot j}{e^{1-\alpha}}
\hspace*{\arraycolsep}=\hspace*{\arraycolsep}
\displaystyle \sum_{0\leq j \leq E}\, \frac{x^{j} \cdot j}{e^{1-\alpha}}.
\end{array} \]

\noindent Subtracting the second line from $\mu$ times the first one,
and then multiplying with $e^{1-\alpha}$ yields a polynomial equation
of order $E$, namely:
\[ \begin{array}{rcl}
0
& = &
\displaystyle \sum_{0\leq j \leq E}\, x^{j} \cdot (i - \mu).
\end{array} \]

\noindent An approximate solution $s>0$ can be found, for instance via
the \pythoninline{nroots} function in \Python's \pythoninline{sympy}
library. We then get a distribution on $\{0,\ldots,E\}$ with mean
$\mu$ and with maximal entropy, as normalisation of:
\begin{equation}
\label{MaxEntropy}
\begin{array}{rcl}
\displaystyle\sum_{0\leq j \leq E} \, s^{j}\,\bigket{j}
& = &
\displaystyle\sum_{0\leq j \leq E} \, e^{-b\cdot j}\,\bigket{j}
   \quad\mbox{ for } b = -\ln(s).
\end{array} 
\end{equation}

\noindent Finally, the formulations~\eqref{DiscExpApprox}
(and~\eqref{MaxEntropy}) lead to a continuous approximation as
exponential with rate $\frac{1}{\mu}$, with probability densitiy
function (pdf) on $\pR$ given by:
\begin{equation}
\label{ContExpApprox}
\begin{array}{rcl}
j 
& \longmapsto &
\frac{1}{\mu} \cdot e^{-\frac{j}{\mu}}
\end{array}
\end{equation}

\begin{figure}
\begin{center}
\hspace*{-1em}\begin{tabular}{cc}
\includegraphics[scale=0.26]{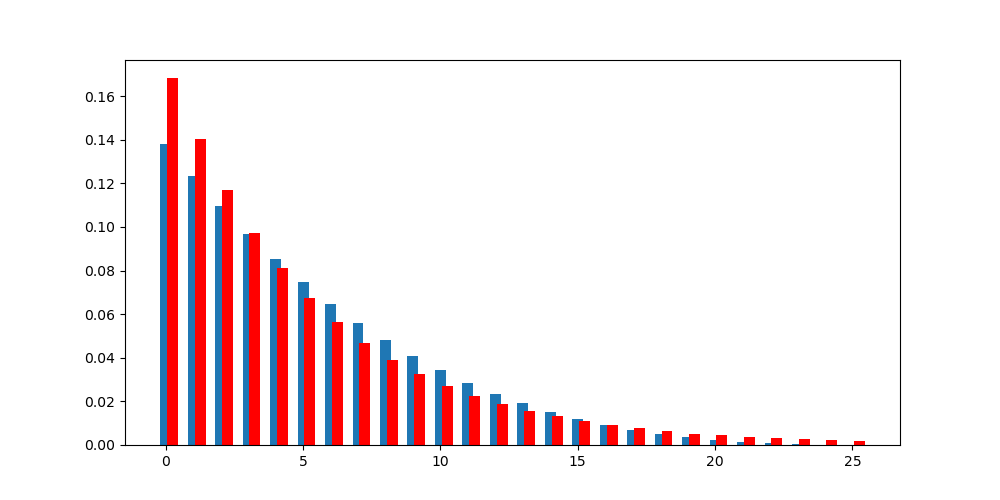}
   \hspace*{-1em} & \hspace*{-1em}
\includegraphics[scale=0.26]{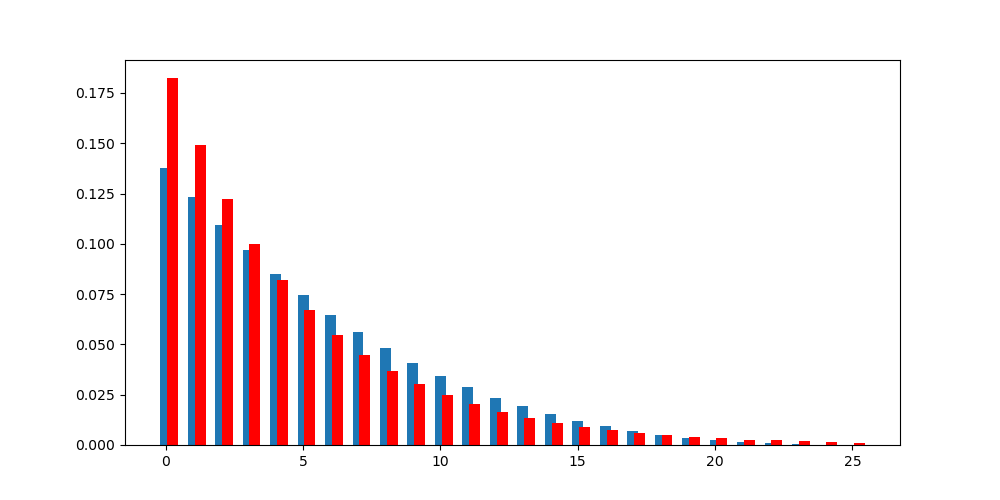}
\\
\includegraphics[scale=0.26]{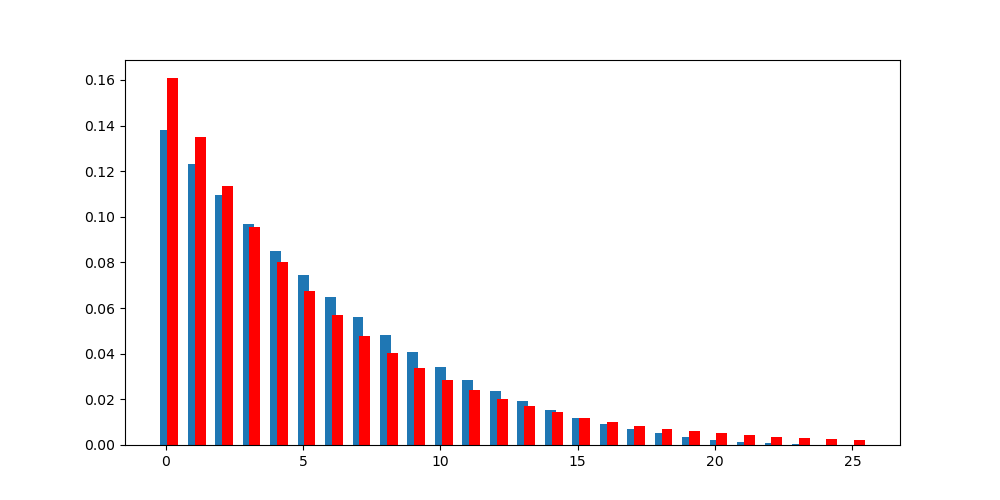}
   \hspace*{-1em} & \hspace*{-1em}
\includegraphics[scale=0.26]{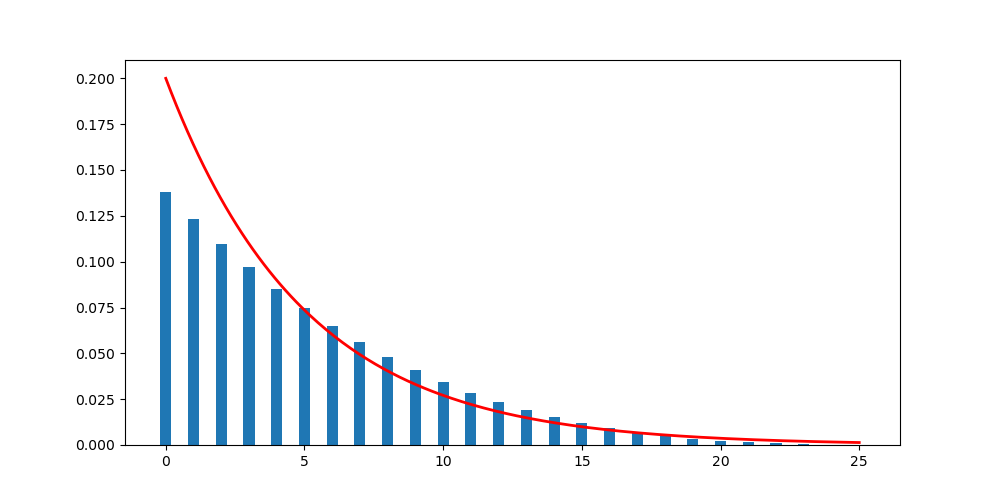}
\end{tabular}
\end{center}
\caption{The Boltzmann-on-energy distribution $\boltzmannene[25](5)$
  with total energy $E=25$ and $K = 5$ particles occurs in all four
  plots as blue bar chart. Its mean is $\mu = \frac{E}{K} = 5$. On the
  top the iterative (normalised) distribution~\eqref{DiscIterApprox}
  and the discrete exponent distribution~\eqref{DiscExpApprox} are
  included in red. Similarly, the bottom contains the maximum entropy
  disribution from~\eqref{MaxEntropy} and the continuous (exponential)
  distribution~\eqref{ContExpApprox}. The maximum entropy distribution
  uses the solution $s = 0.841$, for which $-\ln(s) = 0.173$ is close
  to $\frac{1}{\mu} = 0.2$, which is used in~\eqref{DiscExpApprox}. Its
  entropy is $2.69$, which is bit higher than the entropy $2.67$ of
  the Boltzmann distribution $\boltzmannene[25](5)$.  }
\label{ApproxFig}
\end{figure}

\noindent Figure~\ref{ApproxFig} gives an impression of these four
approximations~\eqref{DiscIterApprox} -- \eqref{ContExpApprox}. When
we look at the KL-divergences between the Boltzmann-on-energy
distribution and the three (discrete)
approximations~\eqref{DiscIterApprox} -- \eqref{MaxEntropy}, then the
last one performs slightly better --- in this example.

\ignore{

E = 25
K = 5
bo = boltzmann_on_energy_bibinom(E)(K)
sp4 = Space([i/4 for i in range(4*E+2)])
mu = bo.expectation()
iter_approx = DState(lambda x: (mu/(mu+1))**x, range_sp(E+1)).flrn()
exp_approx = DState(lambda x: math.exp(-1/mu * x), range_sp(E+1)).flrn()
print ("\nBoltzmann mean: ", mu )
bo4 = Functor(lambda x: x, cod=sp4)(bo)
print(" * Iter approximation mean and distances: ", iter_approx.expectation(), 
      "  ", tvdist(bo, iter_approx), "  ", kldive(bo, iter_approx) )
print(" * Exp approximation mean and distance: ",  exp_approx.expectation(), 
      "  ", tvdist(bo, exp_approx), "  ", kldive(bo, exp_approx) )
sp4 = Space([i/4 for i in range(4*E+2)])
bo4 = Functor(lambda x: x, cod=sp4)(bo)
bo4.dplot1(bar_width=0.4, other_distributions = 
           [ Functor(lambda x: x+0.25, cod=sp4)(iter_approx) ])
bo4.dplot1(bar_width=0.4, other_distributions = 
           [ Functor(lambda x: x+0.25, cod=sp4)(exp_approx) ])
bo.dcplot1(bar_width = 0.4, other_distributions = 
           [ State(lambda x: 1/mu * math.exp(-1/mu * x), R(0,inf)) ])

}

\auxproof{


His value $Nu$ is our energy $e$, where $N$ is our $K$, that is the
number of particles. Thus $u = \frac{e}{N}$ is the mean of the
Boltzmann-on-energy distribution. There is a set $S$ which contains
the sequences that add up to $e$, like in
Lemma~\ref{BoltzmannMltLem}~\eqref{BoltzmannMltLemIm}. A Markov chain
is considered on $S$ where a random element is chosen and added, like
in the UU mode in~\cite{JacobsS25}. The uniform distribution on $S$ is
claimed to be stationary. This yields the multiset-coefficent
distribution when $\acc$ is applied.

The reformulation in
Proposition~\ref{BoltzmannNumAltProp}~\eqref{BoltzmannNumAltPropProj}
is inspired by this approach.

One can put a Markov chain on $S$ by removing and adding single
numbers from single coordinates. Interestingly, this maintains the
sum, whereas this sum preservation is difficult to achieve with the
multiset-based approach from~\cite{JacobsS25}.

In the limit, when intervals get smaller, the space $S$ can be
described via $\nnR$-based multisets, as $\Mlt[i](\finset{N})$.
Interestingly, this is claimed to a cdf, and then a pdf, see the
top of p.3, but I cannnot reconstruct the calculation for the
first equation, although it is intriguing and may give a 
connection between Boltzmann and Dirichlet, e.g. in
experiments:

\tabular{l}
\texttt{a = 0.4}
\\
\texttt{print( D >= Pred(lambda xs: 1 if xs[0] >= a else 0, D.sp) )}
\\
\texttt{print( (1 - a)**2 )}
\end{tabular}


}

\subsection{A Markov chain with its stationary Boltzmann 
   distribution}\label{MarkovChainSubsec}

In~\cite{JacobsS25} several Markov chains are defined on a set
$\natMlt[K](X)$ of multisets of size $K$. Here we describe a more
complicated chain that maintains not only the size of multisets, but
also their sum. Thus, our aim is to define a transition channel
$\shift \colon \natMlt[K,i](\finset{N}) \rightarrow
\Dst\big(\natMlt[K,i](\finset{N})\big)$, forming a Markov chain. The
definition below is a bit technical, but the intuition is simple: we
randomly move one particle one energy level down and then randomly
move one particle one energy level up. Consider the multiset $\varphi
= 1\ket{0} + 2\ket{1} + 3\ket{2} \in \natMlt[6, 8](\finset{3})$ as a
configuration of $6$~particles with a total energy of~$8$. We can turn
it into other configurations, still with the same size~$6$ and
sum~$8$, by carefully downgrading and upgrading single particles.
\begin{itemize}
\item We can downgrade one particle from level~$1$ to level~$0$.  This
  turns $\varphi$ into an intermediate multiset $2\ket{0} + 1\ket{1} +
  3\ket{2}$ via a shift from~$1$ to~$0$. We now have to decide at
  which level to upgrade. This can be at level~$0$ or at level~$1$,
  but not at level~$2$ since that is the highest level. Upgrading at
  level $0$ returns the original multiset $\varphi$, which is
  possible; upgrading at level~$1$ turns the intermediate multiset
  into $2\ket{0} + 4\ket{2}$, via a shift from~$1$ to~$2$. The
  resulting multiset still has size~$6$ and sum~$8$, as required.

\item We can also downgrade one particle in $\varphi$ at level~$2$,
  resulting in an intermediate multiset $1\ket{0} + 3\ket{1} +
  2\ket{2}$, via a shift from~$2$ to~$1$.  Now we can decide to
  upgrade at level~$0$, giving $4\ket{1} + 2\ket{2}$, or at level~$1$,
  giving the original multiset $1\ket{0} + 2\ket{1} + 3\ket{2}$.
\end{itemize}

\noindent As an aside, this downgrading and upgrading works like the
draw-delete and draw-add for multiset partitions in~\cite{Jacobs22c},
in the context of Ewens distributions in population biology.

We are now ready to define the channel $\shift \colon
\natMlt[K,i](\finset{N}) \rightarrow
\Dst\big(\natMlt[K,i](\finset{N})\big)$ in general form. The letters
`$d$' and `$u$' are used to suggest the downgrade and the upgrade, and
$\varphi_{d}$ is what we called the `intermediate' multiset.  Notice
that the downgrade $d$ cannot happen at level~$0$ and the upgrade $u$
not at level~$N\minnetje 1$.
\[ \begin{array}{rcl}
\shift(\varphi)
& \coloneqq &
\displaystyle\frac{\varphi(0)}{\|\varphi\|}\bigket{\varphi} +
   \sum_{0 < d < N, \, \varphi(d) > 0} \; \sum_{0 \leq u < N-1, \, \varphi_{d}(u) > 0} \, 
\\[+1.2em]
& & \hspace*{8em} \displaystyle
   \frac{\varphi(d)}{\|\varphi\|} \cdot 
   \frac{\varphi_{d}(u)}{\|\varphi_{d}\| - \varphi_{d}(N\minnetje 1)} \,
   \Bigket{\varphi_{d} - 1\ket{u} + 1\ket{u\plusje 1}},
\\[+1.0em]
& & \quad\mbox{where } \varphi_{d} \coloneqq 
   \varphi - 1\ket{d} + 1\ket{d\minnetje 1}.
\end{array} \]

\noindent This shift function forms a Markov chain on
the set of multisets $\natMlt[K,i](\finset{N})$ with the crucial
property that the Boltzmann-on-multisets distribution is stationary,
for this chain.  This means that it is a fixed point of the
pushforward~\eqref{PushEqn}, corresponding to a statistical
equilibrium:
\[ \begin{array}{rcl}
\shift_{*}\Big(\boltzmannmlt[N,K](i)\Big)
& = &
\boltzmannmlt[N,K](i).
\end{array} \]

Interestingly, we can turn the shift on multisets into a shift on
numbers, now with the Boltzmann-on-numbers distribution as fixed point.
We use the frequentist learning function as a channel $\flrn
\colon \natMlt[K,i](\finset{N}) \rightarrow \Dst(\finset{N})$ and
construct its `Bayesian inversion' $\flrn^{\dag} \colon \finset{N}
\rightarrow \Dst\big(\natMlt[K,i](\finset{N})\big)$, with the
Boltzmann-on-multisets distribution as prior. We refer
to~\cite{Jacobs21g} for a general account and only provide
the relevant definitions here. We use, for $0\leq j < N$,
\[ \begin{array}{rcl}
\flrn^{\dag}(j)
& \coloneqq &
\displaystyle\sum_{\varphi\in\natMlt[K,i](\finset{N})} \, 
   \frac{\coefm{\varphi} \cdot \varphi(j)}
   {\sum_{\psi\in \natMlt[K,i](\finset{N})} \coefm{\psi}\cdot\psi(j)} 
  \,\bigket{\varphi}.
\end{array} \]

\noindent The shift-on-numbers channel $\finset{N} \rightarrow
\Dst(\finset{N})$ is then obtained via channel composition $\klafter$
from Definition~\ref{PushDef}, as:
\[ \begin{array}{ccrcl}
\flrn \klafter \shift \klafter \flrn^{\dag}
& \quad\mbox{so}\quad &
\Big(\flrn \klafter \shift \klafter \flrn^{\dag}\Big)_{*}
& = &
\flrn_{*} \after \shift_{*} \after \big(\flrn^{\dag}\big)_{*}.
\end{array} \]

\noindent It is not hard to see that it has the Boltzmann-on-numbers
distribution $\boltzmannnum[N,K](i)$ as fixed point (equilibrium),
using that $\big(\flrn^{\dag}\big)_{*}\big(\boltzmannnum[N,K](i)\big)
= \boltzmannmlt[N,K](i)$.

In particular, the Boltzman-on-energy distribution
$\boltzmannene[25](5)$ used in Figure~\ref{ApproxFig} is a stationary
equilibrium, for a certain Markov chain. But it does not exactly match
the maximum entropy distribution~\eqref{MaxEntropy}. Hence, in the
discrete case, there is a mismatch between being an equilibrium (for
the above shift) and have maximal entropy (among distributions with
the same mean).

\subsection{New multivariate distributions from 
   $N$-nomials}\label{NomialDstSubsec}

In Lemma~\ref{NnomialSeqLem}~\eqref{NnomialSeqLemVDM} we have seen a
Vandermonde result for $N$-nomials. Such a Vandermonde result exists
for binomial (as special case) and also for multichoose coefficients,
see~\cite{Jacobs22a}.  They can be generalised to $n$-ary form and are
then best expressed in terms of multisets. We briefly explore this
topic. This requires a bit more notation. For multisets $\varphi, \psi
\in \natMlt(X)$ we write:
\[ \begin{array}{rcl}
\varphi \leq \psi
& \mbox{ when } &
\varphi(x) \leq \psi(x) \mbox{ for each }x\in X
\\
\varphi \leq_{K} \psi
& \mbox{ when } &
\varphi \leq \psi \mbox{ and } \|\varphi\| = K.
\end{array} \]

\noindent Next we extend binomial and multichoose coefficients to
multisets, via pointwise definitions:
\[ \begin{array}{rclcrcl}
\displaystyle\binom{\psi}{\varphi}
& \coloneqq &
\displaystyle\prod_{x\in X} \binom{\psi(x)}{\varphi(x)}
& \qquad\mbox{and}\qquad &
\displaystyle\bibinom{\psi}{\varphi}
& \coloneqq &
\displaystyle\prod_{x\in X} \bibinom{\psi(x)}{\varphi(x)}.
\end{array} \]

\noindent The first definition requires that $\varphi \leq \psi$ and
the second one that $\one \leq \psi$, where $\one = \sum_{x} 1\ket{x}
\in\natMlt(X)$ is the multiset of singletons --- assuming that the set
$X$ is finite.

These definitions are used in~\cite{Jacobs22a,Jacobs21g} for the
following multivariate versions of the Vandermonde properties, for
$\psi\in\natMlt[L](X)$.
\begin{equation}
\label{VDMEqns}
\begin{array}{rclcrcl}
\displaystyle\sum_{\varphi\leq_{K} \psi} \binom{\psi}{\varphi}
& = &
\displaystyle\binom{L}{K}
& \qquad\mbox{and}\qquad &
\displaystyle\sum_{\varphi\in\natMlt[K](X)} \bibinom{\psi}{\varphi}
& = &
\displaystyle\bibinom{L}{K}.
\end{array}
\end{equation}

\noindent These equations are then used to define the multivariate
hypergeometric and P\'olya distributions, as distributions on
multisets of draws, in $\Dst\big(\natMlt[K](X)\big)$.
\[ \begin{array}{rclcrcl}
\hypergeometric[K](\psi)
& \coloneqq &
\displaystyle\sum_{\varphi\leq_{K}\psi} \, 
   \frac{\binom{\psi}{\varphi}}{\binom{L}{K}}\,\bigket{\varphi}
& \quad\mbox{and}\quad &
\polya[K](\psi)
& \coloneqq &
\displaystyle\sum_{\varphi\in\natMlt[K](X)} \, 
   \frac{\big(\binom{\psi}{\varphi}\big)}{\big(\binom{L}{K}\big)}
   \,\bigket{\varphi}.
\end{array} \]

\noindent The hypergeometric distribution captures the draw-and-remove
mode, where a drawn ball is removed from the urn. For this reason the
admissable draws $\varphi$ of size $K$ are required to satisfy
$\varphi \leq_{K} \psi$, so that only balls that are actually in the
urn $\psi$ can be drawn. The P\'olya distribution uses the
draw-and-duplicate mode, where a drawn ball is returned to the urn,
together with an additional ball of the same colour,
see~\cite{Jacobs22a} for details.  These two distributions both
interact well with the pushforward~\eqref{PushEqn} of frequentist
learning, in the sense that:
\begin{equation}
\label{FlrnDrawEqn}
\begin{array}{rclcrcl}
\flrn_{*}\Big(\hypergeometric[K](\psi)\Big)
& = &
\flrn(\psi)
& \quad\mbox{and}\quad &
\flrn_{*}\Big(\polya[K](\psi)\Big)
& = &
\flrn(\psi).
\end{array}
\end{equation}

\noindent These are fundamental results, see also~\cite{Jacobs21g},
that express that applying frequentist learning `in probability' to
the draws from urn $\psi$ yields the distribution associated with the
original urn, again via frequentist learning.

These observations prepare us for a similar situation for
$N$-nomials. We can generalise the Vandermonde result in
Lemma~\ref{NnomialSeqLem}~\eqref{NnomialSeqLemVDM} to multivariate
form via multisets. Let a number $N\geq 1$ and a multiset
$\psi\in\natMlt[K](X)$ of size $K$ be given. We write $(N\minnetje
1)\cdot \psi$ for the $N\minnetje 1$ fold sum $\psi + \cdots + \psi$
of multisets, so that $(N\minnetje 1)\cdot \psi = \sum_{x} (N\minnetje
1)\cdot \psi(x) \bigket{x}$. For a number $0\leq i \leq (N\minnetje
1)\cdot K$ and a multiset $\varphi\in\natMlt[i](X)$ with $\varphi \leq
(N\minnetje 1)\cdot \psi$, we define a product of $N$-nomials:
\begin{equation}
\label{MltNomialEqn}
\begin{array}{rcl}
C_{N}\big(\psi, \varphi)
& \coloneqq &
\displaystyle\prod_{x\in X} C_{N}\big(\psi(x), \varphi(x)\big).
\end{array}
\end{equation}

\noindent Then one can show, analogously to~\eqref{VDMEqns},
\[ \begin{array}{rcl}
\displaystyle\sum_{\varphi\leq_{i} (N-1)\cdot \psi} \,
   C_{N}\big(\psi, \varphi\big)
& = &
C_{N}\big(K, i\big).
\end{array} \]

\ignore{

# Multivariate Vandermonde for N-nomials
# number of energy levels
M = 4
X = alphaspace(M)
L = 10
sizes = random_multiset(L-M)(X) + unit_multiset(X)
N = random.randint(2,6)
E = random.randint(0, (N-1)*L)

print("\nNomial for N =", N, "and E =", E, "and sizes: ", sizes)
print( nomial(N, L, E) )
print( sum([ list_multiplication([ nomial(N, sizes(x), e(x)) 
                                   for x in X.iter_discrete() ])
             for e in Multisets(E)(X) if e <= (N-1) * sizes ]) )

}

\noindent This equation now gives rise to a new multivariate
distribution of nomials, namely:
\begin{equation}
\label{NomialDstEqn}
\begin{array}{rcl}
\nomial[i](\psi)
& \coloneqq &
\displaystyle\sum_{\varphi\leq_{i} (N-1)\cdot \psi} \,
   \frac{C_{N}(\psi, \varphi)}{C_{N}(L, i)} \, \bigket{\varphi}
\qquad\mbox{where } N = \setsize{X}.
\end{array} 
\end{equation}

\noindent This distribution is well-behaved, in the sense that
it satisfies the analogue of~\eqref{FlrnDrawEqn}:
\[ \begin{array}{rcl}
\flrn_{*}\Big(\nomial[i](\psi)\Big)
& = &
\flrn(\psi).
\end{array} \]

\noindent It is unclear if there is an urn model for this nomial
distribution in~\eqref{NomialDstEqn}. It is some form of
generalisation of the hypergeometric distribution, since in the binary
case, when the underlying set $X$ has two elements, one has
$\nomial[i](\psi) = \hypergeometric[i](\psi)$. But what it does more
generally is an open question. We conclude with an illustratation,
with urn / multiset $\psi = 1\ket{a} + 5\ket{b} + 3\ket{c}$; then:
\[ \begin{array}{rcl}
\lefteqn{\nomial[15](\psi)}
\\
& = &
\frac{7}{156}\Bigket{2\ket{a} + 10\ket{b} + 3\ket{c}} + 
\frac{5}{26}\Bigket{2\ket{a} + 9\ket{b} + 4\ket{c}} + 
\frac{1}{26}\Bigket{1\ket{a} + 10\ket{b} + 4\ket{c}}
\\[+0.2em]
& & \quad +\;
\frac{15}{52}\Bigket{2\ket{a} + 8\ket{b} + 5\ket{c}} + 
\frac{5}{52}\Bigket{1\ket{a} + 9\ket{b} + 5\ket{c}} + 
\frac{1}{52}\Bigket{10\ket{b} + 5\ket{c}}
\\[+0.2em]
& & \quad +\;
\frac{5}{26}\Bigket{2\ket{a} + 7\ket{b} + 6\ket{c}} + 
\frac{5}{52}\Bigket{1\ket{a} + 8\ket{b} + 6\ket{c}} + 
\frac{5}{156}\Bigket{9\ket{b} + 6\ket{c}}.
\end{array} \]

\noindent This is neither a hypergeometric nor a P\'olya distribution.

The $N$-nomials for multisets~\eqref{MltNomialEqn} can be used to
define a multivariate version of the Boltzmann-on-multisets
distribution, giving a `Boltzmann-on-multiple-multisets' distribution.
This starts with $N$ energy levels and a total energy $i$, and with a
multiset $\psi\in\natMlt[K](X)$ of sizes, where the finite set $X =
\{x_{1}, \ldots, x_{L}\}$ describes the different kinds $x_j$ of
particles at hand.
\[ \begin{array}{rcl}
\hspace*{-0.2em}\boltzmannmltmlt[N, \psi](i)
& \!\coloneqq\! &
\displaystyle\hspace*{-0.4em}\sum_{\begin{array}{c}
    \\[-2.5em]
   \scriptstyle\varphi_{1} \in \natMlt[\psi(x_{1})](\finset{N}), \ldots, 
   \varphi_{L} \in \natMlt[\psi(x_{L})](\finset{N}) 
   \\[-0.7em]
   \scriptstyle \som(\varphi_{1}) + \cdots + \som(\varphi_{L}) = E
  \end{array}} \!
  \frac{\coefm{\varphi_1} \cdot \ldots \cdot \coefm{\varphi_{L}}}
     {C_{N}(K, i)} \, \Bigket{\varphi_{1}, \ldots, \varphi_{L}}.
\end{array} \]

\noindent This is a distribution on $\natMlt[\psi(x_{1})](\finset{N})
\times \cdots \times \natMlt[\psi(x_{L})](\finset{N}$. It can be
turned into a distribution on $L$-tuples of numbers via multiple
applications of frequentist learning. It is unclear whether such
multivariate Boltzmann distributions are useful in physics or in other
areas.

\auxproof{
We check that this yields an actual distribution:
\[ \begin{array}{rcl}
\lefteqn{\sum_{\begin{array}{c}
    \\[-2.5em]
   \scriptstyle\varphi_{1} \in \natMlt[\psi(x_{1})](\finset{N}), \ldots, 
   \varphi_{L} \in \natMlt[\psi(x_{L})](\finset{N}) 
   \\[-0.7em]
   \scriptstyle \som(\varphi_{1}) + \cdots + \som(\varphi_{L}) = i
  \end{array}} \!
  \frac{\coefm{\varphi_1} \cdot \ldots \cdot \coefm{\varphi_{L}}}
     {C_{N}(K, i)}}
\\
& = &
\displaystyle \sum_{i_{1}, \ldots, i_{L}, \, i_{1} + \cdots + i_{L} = i} \;
   \sum_{\begin{array}{c}
    \\[-2.5em]
   \scriptstyle\varphi_{1} \in \natMlt[\psi(x_{1})](\finset{N}), \ldots, 
   \varphi_{L} \in \natMlt[\psi(x_{L})](\finset{N}) 
   \\[-0.7em]
   \scriptstyle \som(\varphi_{1}) = i_{1}, \ldots, \som(\varphi_{L}) = i_{L}
  \end{array}} \!
  \frac{\coefm{\varphi_1} \cdot \ldots \cdot \coefm{\varphi_{L}}}
     {C_{N}(K, i)}
\\[+1.4em]
& = &
\displaystyle\sum_{\chi \leq_{i} (N-1)\cdot \psi} \,
   \sum_{\begin{array}{c}
    \\[-2.5em]
   \scriptstyle\varphi_{1} \in \natMlt[\psi(x_{1})](\finset{N}), \ldots, 
   \varphi_{L} \in \natMlt[\psi(x_{L})](\finset{N}) 
   \\[-0.7em]
   \scriptstyle \som(\varphi_{1}) = \chi(x_{1}), \ldots, 
   \som(\varphi_{L}) = \chi(x_{L})
  \end{array}} \!
  \frac{\coefm{\varphi_1} \cdot \ldots \cdot \coefm{\varphi_{L}}}
     {C_{N}(K, i)}
\\[+1.4em]
& = &
\displaystyle\sum_{\chi \leq_{i} (N-1)\cdot \psi} \,
   \frac{C_{N}(\psi(x_{1}), \chi(x_{1})) \cdot \ldots \cdot 
   C_{N}(\psi(x_{L}), \chi(x_{L}))}{C_{N}(K, i)}
\\[+1.4em]
& = &
\displaystyle\sum_{\chi \leq_{i} (N-1)\cdot \psi} \,
   \frac{C_{N}(\psi, \chi)}{C_{N}(K, i)}
\\
& = &
1.
\end{array} \]
}

\section{Conclusions}\label{ConclusionSec}

In the area of statistical mechanics one finds intuitive concrete
descriptions of particle configurations, with a fixed number of
particles and with a fixed total energey sum, see
\textit{e.g.}~\cite{DillB10}. This paper elaborates the underlying
combinatorics and statistics in a new manner, via $N$-nomials and
multisets. Physicists often skip this combinatorial part, or are even
reluctant to consider it, as in~\cite{Ramshaw18}. Our approach leads
to general formulas that capture the concrete descriptions in the
physics literature as instantiations. We dare to think that the
systematic use of multisets offers an enrichment not only of
combinatorics but also of statistical mechanics.


\subsection*{Acknowledgements}

\noindent Thanks to M\'ark Sz\'eles and Dario Stein for their helpful
comments on an earlier version of this paper.

\appendix

\section{Addendum to Example~\ref{PhysicsMltEx}}

We elaborate the illustration that was briefly mentioned at the end of
Example~\ref{PhysicsMltEx}, coming
from~\cite{TiplerLResources}. There, a table ``BD-1'' is described
with ``States and occupation probabilities for six particles with
total energy $8\Delta E$''. This $8\Delta E$ correspond in our setting
to the situation in Definition~\ref{BoltzmannEneDef} with $E = 8$ and
$N = E+1 = 9$, in presence of $K=6$ particles.  The table ``BD-1''
lists the 20 multisets $\varphi\in\natMlt[6](\finset{9})$, called
macrostates, with $\som(\varphi) = 8$, together with their multiset
coefficients $\coefm{\varphi}$, called ``number of microstates'', that
is, number of sequences that accumulate to $\varphi$, see
Lemma~\ref{MltLem}~\eqref{MltLemAcc}. The sum of these 20 multiset
coefficients $\coefm{\varphi}$ is listed as 1287, which in our setting
appears as $C_{9}(6,8) = 1287$, using~\eqref{NnomialMltEqn}.

The resulting Boltzmann-on-multisets distribution is big. With the
order of multisets / macrostates as in~\cite[Table
  BD-1]{TiplerLResources}, it is:
\[ \begin{array}{rcl}
\lefteqn{\boltzmannmlt[9,6](8)
\hspace*{\arraycolsep}\smash{\stackrel{\eqref{BoltzmannMltEqn}}{=}}\hspace*{\arraycolsep}
\displaystyle\sum_{\varphi\in\natMlt[6](\finset{9}), \, \som(\varphi) = 8} \, 
   \frac{\coefm{\varphi}}{C_{9}(6,8)} \, \bigket{\varphi}}
\\[+0.4em]
& = &
\frac{2}{429}\Bigket{5\ket{0} + 1\ket{8}} +
\frac{10}{429}\Bigket{4\ket{0} + 1\ket{1} + 1\ket{7}} + 
\frac{10}{429}\Bigket{4\ket{0} + 1\ket{2} + 1\ket{6}}
\\[+0.2em]
& & \quad +\,
\frac{10}{429}\Bigket{4\ket{0} + 1\ket{3} + 1\ket{5}} + 
\frac{5}{429}\Bigket{4\ket{0} + 2\ket{4}} + 
\frac{20}{429}\Bigket{3\ket{0} + 2\ket{1} + 1\ket{6}} 
\\[+0.2em]
& & \quad +\,
\frac{20}{429}\Bigket{3\ket{0} + 2\ket{2} + 1\ket{4}} + 
\frac{20}{429}\Bigket{3\ket{0} + 1\ket{2} + 2\ket{3}}
\\[+0.2em]
& & \quad +\,
\frac{40}{429}\Bigket{3\ket{0} + 1\ket{1} + 1\ket{2} + 1\ket{5}} +
\frac{40}{429}\Bigket{3\ket{0} + 1\ket{1} + 1\ket{3} + 1\ket{4}}
\\[+0.2em]
& & \quad +\,
\frac{5}{429}\Bigket{2\ket{0} + 4\ket{2}} + 
\frac{10}{143}\Bigket{2\ket{0} + 2\ket{1} + 2\ket{3}}
\\[+0.2em]
& & \quad +\,
\frac{20}{143}\Bigket{2\ket{0} + 1\ket{1} + 2\ket{2} + 1\ket{3}} + 
\frac{20}{143}\Bigket{2\ket{0} + 2\ket{1} + 1\ket{2} + 1\ket{4}}
\\[+0.2em]
& & \quad +\,
\frac{20}{429}\Bigket{2\ket{0} + 3\ket{1} + 1\ket{5}} + 
\frac{10}{429}\Bigket{1\ket{0} + 4\ket{1} + 1\ket{4}}
\\[+0.2em]
& & \quad +\,
\frac{40}{429}\Bigket{1\ket{0} + 3\ket{1} + 1\ket{2} + 1\ket{3}} + 
\frac{20}{429}\Bigket{1\ket{0} + 2\ket{1} + 3\ket{2}}
\\[+0.2em]
& & \quad +\,
\frac{5}{429}\Bigket{4\ket{1} + 2\ket{2}} + 
\frac{2}{429}\Bigket{5\ket{1} + 1\ket{3}}.
\end{array} \]


At the end of~\cite[Table BD-1]{TiplerLResources} there is a sequence
of 9~numbers:
\[ 2.31 \quad 1.54 \quad 0.98 \quad 0.59 \quad 0.33 \quad 0.16 
   \quad 0.07 \quad 0.02 \quad 0.005 \]

\noindent We can reconstruct this sequence via the Boltzmann-on-energy
distribution. First:
\[ \begin{array}{rcl}
\lefteqn{\boltzmannene[8](6)
\hspace*{\arraycolsep}\smash{\stackrel{\eqref{BoltzmannEneEqn}}{=}}\hspace*{\arraycolsep}
\boltzmannnum[9,6](8)
\hspace*{\arraycolsep}\smash{\stackrel{\eqref{BoltzmannEneEqn}}{=}}\hspace*{\arraycolsep}
\displaystyle\sum_{0\leq j \leq 8} \, 
      \frac{\big(\binom{5}{8-j}\big)}{\big(\binom{6}{8}\big)}\,\bigket{j}}
\\[+0.4em]
& = &
\frac{5}{13}\bigket{0} + 
\frac{10}{39}\bigket{1} + 
\frac{70}{429}\bigket{2} + 
\frac{14}{143}\bigket{3} + 
\frac{70}{1287}\bigket{4}
\\[+0.2em]
& & \quad +\,
\frac{35}{1287}\bigket{5} + 
\frac{5}{429}\bigket{6} + 
\frac{5}{1287}\bigket{7} + 
\frac{1}{1287}\bigket{8}.
\end{array} \]


\noindent If we take the un-normalised version, times $K=6$, and write
the probabilities as decimals we recognise the above sequence as
probabilities:
\[ \begin{array}{rcl}
6\cdot \boltzmannene[8](6)
& \approx &
2.31\bigket{0} + 
1.54\bigket{1} + 
0.979\bigket{2} + 
0.587\bigket{3}
\\[+0.2em]
& & \quad +\,
0.326\bigket{4} + 
0.163\bigket{5} + 
0.0699\bigket{6} + 
0.0233\bigket{7} + 
0.00466\bigket{8}.
\end{array} \]

\noindent We see that this paper provides a general framework in which
such illustrations fit.


\ignore{

print("\n9-nomial: ", nomial(9, 6, 8) )
bom = boltzmann_on_multisets(9, 6, frac=True)(8)
#print("\nNumber of multisets / configurations: ", len(bom.support_list()) )
print("")
print( bom )

print("\nBoltzmann-on-numbers, also Boltzmann-on-energy")
bon = boltzmann_on_numbers_rec(9, 6, frac=True)(8)
boe = boltzmann_on_energy_bibinom(8, frac=True)(6)
print("")
print( bon )
print("")
print( boe )

print("\nUnnormalised, times K=6")
print( 6 * boe.no_frac() )

}

\auxproof{
At the end of the supplementary note two other outcomes are
computed.
\begin{itemize}
\item The averages per energy level, without taking any probabilities
into account, are obtained as:
\[ \begin{array}{rcl}
\displaystyle\sum_{\varphi\in\supp(\boltzmannmlt[9,6](8))} \,
   \frac{\varphi}{\setsize{\supp(\boltzmannmlt[9,6](8))}}
& = &
\displaystyle\sum_{\varphi\in\supp(\boltzmannmlt[9,6](8))} \,
   \frac{\varphi}{20}.
\end{array} \]

\noindent This is associated with Bose-Einstein.

\item Then there is also the Pauli exclusion principle, which means
  that only those $\varphi$ are admissable with $\varphi(j) \leq 2$
  for each $j$. This holds only for the multisets occurring as 12th,
  13th, 14th in the above description of $\boltzmannmlt[9,6](8)$,
that is, for:
\[ 2\ket{0} + 2\ket{1} + 2\ket{3}
\quad
2\ket{0} + 1\ket{1} + 2\ket{2} + 1\ket{3}
\quad
2\ket{0} + 2\ket{1} + 1\ket{2} + 1\ket{4} \]

\noindent Again, the average of this is computed as in the previous
bullet.
\end{itemize}
}



\end{document}